\documentclass[11pt]{amsart}
\usepackage{amssymb,amsfonts,amsmath}
\usepackage[dvips]{graphicx}
\newtheorem{theo}{Theorem}[section]
\newtheorem{prop}[theo]{Proposition}
\newtheorem{lemma}[theo]{Lemma}

\newtheorem{conj}[theo]{Conjecture}

\newcommand{\sm}{\setminus}

\newcommand{\eps}{\varepsilon}

\newcommand{\G}{{\mathcal G}}
\newcommand{\F}{{\mathcal F}}
\newcommand{\T}{{\mathcal T}}

\newcommand{\Q}{{\mathcal Q}}
\def\proof{\noindent{\bf Proof.}\ }

\def\COMMENT#1{}
\def\TASK#1{}
\def\noproof{{\unskip\nobreak\hfill\penalty50\hskip2em\hbox{}\nobreak\hfill%
       $\square$\parfillskip=0pt\finalhyphendemerits=0\par}\goodbreak}
\def\endproof{\noproof\bigskip}
\newdimen\margin   
\def\textno#1&#2\par{%
   \margin=\hsize
   \advance\margin by -4\parindent
          \setbox1=\hbox{\sl#1}%
   \ifdim\wd1 < \margin
      $$\box1\eqno#2$$%
   \else
      \bigbreak
      \hbox to \hsize{\indent$\vcenter{\advance\hsize by -3\parindent
      \sl\noindent#1}\hfil#2$}%
      \bigbreak
   \fi}
\def\proof{\removelastskip\penalty55\medskip\noindent{\bf Proof. }}
\def\enddiscard{}
\long\def\discard#1\enddiscard{}

\title{An approximate version of Sumner's universal tournament conjecture}
\author{Daniela K\"uhn \and Richard Mycroft \and Deryk Osthus}
\date{}
\thanks {D.~K\"uhn was partially supported by the EPSRC, grant
no.~EP/F008406/1. D.~Osthus was partially supported by the EPSRC, grant
no.~EP/E02162X/1 and~EP/F008406/1.}

\begin{document}

\vspace*{-0.8cm}
\begin{abstract}Sumner's universal tournament conjecture states that any tournament on $2n-2$
vertices contains a copy of any directed tree on $n$ vertices. We prove an
asymptotic version of this conjecture, namely that any tournament on $(2+o(1))n$
vertices contains a copy of any directed tree on~$n$ vertices. In addition, we
prove an asymptotically best possible result for trees of bounded degree, namely
that for any fixed $\Delta$, any tournament on $(1+o(1))n$ vertices contains a
copy of any directed tree on $n$ vertices with maximum degree at most $\Delta$.
\end{abstract}

\maketitle

\vspace*{-0.6cm}
\section{Introduction}\label{sec:preliminaries}

\subsection{Introduction}

A tournament is an orientation of a complete graph.
One of the most well-known problems on tournaments is Sumner's universal
tournament conjecture,
which was posed in 1971 (see e.g.~\cite{RW,W}).
\begin{conj} \label{main_conjecture}
Let $T$ be a directed tree on $n$ vertices. Then every tournament on $2n-2$
vertices
contains a copy of $T$.
\end{conj}
The following simple example shows that the bound would be best possible:
let $G$ be a regular tournament on $2n-3$
vertices (so every vertex has $n-2$ outneighbours), and let $T$ be a star with all
edges directed outwards. Then the central vertex of $T$ has $n-1$ outneighbours,
and so $G$ does not contain a copy of $T$.

A large number of partial results towards Sumner's conjecture have been obtained.
Let $f(n)$ denote the
smallest integer such that any tournament on $f(n)$ vertices contains any
directed tree on $n$ vertices. So Conjecture~\ref{main_conjecture} states that
$f(n) = 2n-2$. Chung (see \cite{W}) observed that $f(n) \leq n^{1+o(1)}$, and
Wormald~\cite{W} improved this bound to $f(n) \leq n \log_2 (2n/e)$. The first
linear bound on $f(n)$ was established by H\"aggkvist and Thomason~\cite{HT},
who showed that $f(n) \leq 12n$, and also that $f(n) \leq (4+o(1))n$.
Havet~\cite{H} showed that $f(n) \leq 38n/5-6$, and then Havet and Thomass\'e~\cite{HTh}
used the notion of median orders to improve this to $f(n) \leq (7n-5)/2$.
The current best bound is due to El Sahili~\cite{ES}.

\begin{theo} [\cite{ES}] \label{best_so_far}
Let $T$ be a directed tree on $n$ vertices. Then every tournament on $3n-3$
vertices
contains a copy of $T$.
\end{theo}

The conjecture has also been verified for some classes of trees, such as
directed paths. Indeed, a classical result of Redei~\cite{redei} implies that we
 can even find a spanning directed path in any
tournament.

\begin{theo} [\cite{redei}] \label{dir_path_embed}
For any positive integer $n$, any tournament on $n$ vertices contains a directed
path on $n$ vertices.
\end{theo}
Thomason~\cite{T} proved a much stronger result, namely that whenever $n$ 
is sufficiently large, every tournament on $n$ vertices contains every orientation 
of the path on $n$ vertices (this was a conjecture of Rosenfeld). Havet and 
Thomass\'e~\cite{HTrosenfeld} showed that this even holds for all $n \neq {3,5,7}$. 
Reid and Wormald~\cite{RW} also proved Sumner's conjecture for other (very restricted)
classes of trees.
Havet and Thomass\'e~\cite{HTh} proved that Conjecture~\ref{main_conjecture} holds for
arborescences, i.e.~where~$T$ has a specified root $r$ so that either every edge
of~$T$ is
directed towards $r$, or every edge of~$T$ is directed away from $r$.

We will prove an approximate version of Sumner's conjecture. We
also prove an asymptotically sharp bound for trees with bounded maximum
degree.
\begin{theo} \label{main}
Let $\alpha >0$. Then the following properties hold.
\begin{enumerate}
\item There exists $n_0$ such that for any $n \geq n_0$, any
tournament $G$ on $2(1+\alpha)n$ vertices contains any directed tree $T$ on $n$ vertices.
\item Let $\Delta$ be any positive integer. Then there exists $n_0$ such that
for any $n \geq n_0$, any tournament $G$ on $(1+\alpha)n$
vertices contains any directed tree $T$ on $n$ vertices with $\Delta(T) \leq \Delta$.
\end{enumerate}
\end{theo}
In~\cite{exactsumner}, we prove Sumner's conjecture for large~$n$.
The proof relies on the results (and not just the methods)
that we prove in this paper.

Part~(2) of Theorem~\ref{main} implies that
Sumner's conjecture is true with room to spare for large trees of small maximum
degree. The following example shows that (2) is best possible in the sense
that the term $\alpha n$ cannot be completely omitted:
take a regular tournament~$H_1$ on $2k-1$ vertices,
take an arbitrary tournament $H_2$ on $n-k-1$ vertices and obtain a tournament~$G$
on $n+k-2$ vertices from $H_1\cup H_2$ by adding all edges directed from $H_1$ to~$H_2$.
Also, let~$T$ be the tree on $n$ vertices obtained from a directed path on
$n-k$ vertices by adding~$k$ extra vertices which all send an edge to the initial vertex of the
path. Then~$G$ contains no copy of~$T$.
(We are grateful to P.~Allen and O.~Cooley for pointing out this example to us.)
It would be interesting to know whether the term $\alpha n$ can be reduced to
a constant depending only on~$\Delta$.

Another class of trees where Sumner's conjecture can  be strengthened are
trees with few leaves. The first result in this direction was proved by
H\"aggkvist and Thomason~\cite{HT}.
Havet and Thomass\'e (see~\cite{Havetconj}) then proposed
the following generalization of Sumner's conjecture.
\begin{conj} [\cite{Havetconj}] \label{HTleaves}
Let $T$ be a directed tree on $n$ vertices with $k$ leaves. Then every tournament on $n+k-1$
vertices contains a copy of $T$.
\end{conj}
C\'eroi and Havet~\cite{Havet3leaves}
proved that this conjecture holds for $k \le 3$, from which they deduced that
Sumner's conjecture
holds for all trees with at most 4 leaves.

For our proof of Theorem~\ref{main} we introduce a decomposition of an
arbitrary tournament which searches for dense expanding
subgraphs. We then introduce a randomized algorithm
for embedding arbitrary trees into such dense expanding graphs.
Both tools may be useful for other problems. For example,
it would be interesting to know whether our methods can be extended to prove an
approximate version
of Conjecture~\ref{HTleaves}.

\subsection{Outline of the proof.}

The notion of a robust outexpander (which was introduced for dense graphs in~\cite{KOT}) is
crucial to the proof. Informally, a digraph $G$ is a robust
outexpander if for any set $S \subseteq V(G)$ which is not too large or too
small, the number of vertices with many
inneighbours in $S$ is substantially bigger than $|S|$. K\"uhn,
Osthus and Treglown~\cite{KOT} showed that any robust outexpander $G$ of
linear minimum semidegree contains a Hamilton cycle. (Here the minimum semidegree 
is the minimum of the minimum indegree and the minimum
outdegree.) Applying this to the `reduced digraph'
obtained from the Szemer\'edi regularity lemma, this implies that we can split
most of the vertices of $G$ into sets $V_1,
V_2, \dots, V_k$ so that the set of edges from $V_i$ to $V_{i+1}$ for each~$i$
(addition of the indices taken modulo $k$) forms a quasirandom and dense
bipartite graph. As we shall see, this structure is very useful for embedding
trees. On the other hand, it is
easy to show that if a tournament~$G$ is not a robust outexpander of linear
minimum semidegree, then the vertices of $G$ can be split into two parts so that
almost all of the edges between the two parts are directed the same way (see
Lemma~\ref{not_robust_expander_split}). We
shall then consider whether either of these two parts are robust outexpanders,
and so on.

To begin, in Section~\ref{subsec:defs} we shall define the concepts we shall
use, and prove various lemmas which will be of use to us later on. Then in
Sections~\ref{sec:bounded_robust_exp} and~\ref{sec:unbounded_robust_exp} we show
that Theorem~\ref{main} holds with the added condition that $G$ is a robust
outexpander of linear minimum semidegree. Indeed, in
Section~\ref{sec:bounded_robust_exp}, we consider the case where the
tournament~$G$ is a robust outexpander of linear minimum semidegree
on $(1+\alpha)n$
vertices, and $T$ is a directed tree on $n$ vertices of bounded maximum degree.
As described above, we can split most of the vertices of $G$ into clusters
$V_1, V_2, \dots, V_k$ so that the set of edges from $V_i$ to $V_{i+1}$ is
quasirandom and dense for each $i$. Given this structure on $G$, one attempt
to embed $T$ in $G$ would be to embed each vertex $t \in T$ in the cluster
either preceding or succeeding the cluster containing the parent $t'$ of $t$,
according to the direction of the edge between $t$ and $t'$. However, for many
trees this method will fail to give an approximately uniform allocation of
vertices of $T$ to the clusters of $G$, which we require for the embedding to be
successful. Instead, we modify this method so that each vertex is embedded as
above with probability 1/2 and is embedded in the same cluster as its parent
with probability 1/2. We show that with
high probability this randomised algorithm will indeed give an approximately
uniform allocation of vertices of $T$ to the clusters of $G$, and so will
successfully embed $T$ in $G$.

In Section~\ref{sec:unbounded_robust_exp} we begin by strengthening the result
from Section~\ref{sec:bounded_robust_exp}, showing that if $T$ is a directed
tree on $n$ vertices of
bounded maximum degree, and $G$ is a tournament on $(1+\alpha)n$ vertices whose
reduced graph defined on the clusters $V_1,\dots, V_k$ contains a Hamilton cycle,
then we can embed $T$ in $G$ so that
the vertices of a chosen small set $H \subseteq V(T)$ are embedded within a
specified set
$U \subseteq V(G)$. To do this, we embed all
vertices `far' from $H$ by the method described above, which ensures that the
vertices of $T$ are allocated approximately uniformly amongst the
clusters of~$G$.
The remaining vertices of $T$ are instead embedded to ensure that every
vertex of $H$ is embedded within $U$. This result allows us to consider
directed trees $T$ of unbounded
maximum degree. Indeed, we define for a tree~$T$ a `core tree'~$T_c$, which has
the properties that $T_c$ has bounded maximum degree, but each component of
$T-T_c$ is small. This enables us to show that any tournament~$G$ which is a
robust outexpander of linear minimum semidegree on $(2+\alpha)n$ vertices contains any
directed tree on $n$ vertices. To do this, we again split most of the vertices of $G$
into sets $V_1, V_2, \dots, V_k$ as described above. We then choose subsets
$V_i' \subseteq V_i$ at random so that $|\bigcup_i V'_i|$ is roughly equal to $|T_c|$,
and embed $T_c$ into these subsets (actually we first extend $T_c$ to an 
`extended tree' $T_{\textrm{ext}}$ and embed $T_{\textrm{ext}}$ into these subsets), using the
strengthened result for bounded degree trees to restrict certain vertices of
$T_c$ to vertices of~$G$ with many inneighbours and outneighbours in $G -
\bigcup_i V_i'$. Since each component of $T-T_c$ is small, this will allow us to embed
the components of $T-T_c$ one by one in the unoccupied vertices of $G$ to
complete the embedding of $T$ in $G$.

It is a simple exercise to demonstrate that any transitive tournament on $n$
vertices contains any directed tree on $n$ vertices. In
Section~\ref{sec:transitive}, we prove an analogue of this for
almost-transitive tournaments~$G$. This means that the vertices of $G$ can be ordered so
that almost all of the edges of $G$ are directed towards the endvertex
which is greater in this order. We show that if $G$ is an almost-transitive
tournament on $(1+\alpha)n$ vertices and $T$ is a directed tree on $n$ vertices
then $G$ contains~$T$.

Finally, in Section~\ref{proof}, we shall use the robust outexpander dichotomy
to prove Theorem~\ref{main}. Here we shall describe the proof of the first
statement; the proof of the second is very similar. So let $G$ be a tournament
on $2(1+\alpha)n$ vertices and let $T$ be a directed tree on $n$ vertices.
If $G$ is
a robust outexpander of linear minimum semidegree, then our results of
Sections~\ref{sec:bounded_robust_exp} and~\ref{sec:unbounded_robust_exp} show
that $G$ contains $T$, as desired. On the other hand, if $G$ is not a robust
outexpander of linear minimum semidegree then we may split $G$ into two parts as
described above. We now
examine the larger of these two parts. If this is a robust outexpander of linear
minimum semidegree then we stop; otherwise we again split this part into two. We
continue in this fashion, always choosing the largest part of $G$, stopping if
this is a robust outexpander and splitting it into two smaller parts if not.
If we continue this process but do not find a robust outexpander of linear
minimum semidegree, then $G$ must be almost transitive. Indeed, each time we split $G$
most of the edges across the split are directed the same way. So once all of the
parts of $G$ are sufficiently small, we can be sure that for some ordering of
the vertices of $G$, almost all of the edges of~$G$ are directed according to
this order. So by the result from Section~\ref{sec:transitive}, $G$ contains $T$, as desired.

So suppose instead that at some stage we stop because the largest part of $G$ is
a robust outexpander of linear minimum semidegree. Then we divide $T$ into parts
to be embedded amongst the
parts of $G$, so that each part of $G$ receives a part of $T$ approximately proportional
to its size. However, the robust outexpander part of $G$ will actually receive
slightly more vertices of $T$ than it would from a proportional split. The results from
Sections~\ref{sec:bounded_robust_exp} and~\ref{sec:unbounded_robust_exp} guarantee that
this part of $T$ can still be embedded into the corresponding part of $G$. Since
then the other parts of $G$ will receive slightly fewer vertices of $T$ than they
would from a proportional split it will be possible to embed the remainder of~$T$.

\section{Definitions.} \label{subsec:defs}

\subsection{Notation}

For a graph $G$, we shall write $V(G)$ for the vertex set of $G$, and $|G|$ for
the number of vertices of $G$. $E(G)$ denotes the set of edges of $G$, and
$e(G) := |E(G)|$. Similarly for sets $X, Y \subseteq V(G)$, $e(X, Y)$ denotes the
number of edges between $X$ and $Y$. We shall sometimes write $v \in G$ to mean
$v \in V(G)$. The
\emph{degree} of a vertex $v \in G$, denoted $d(v)$, is the number of edges $e
\in E(G)$ incident
to $v$. We denote the minimum and maximum degree (taken over all vertices of
$G$) by $\delta(G)$, and $\Delta(G)$ respectively. The \emph{distance} $d(u, v)$
between vertices $u, v \in G$ is the length of the shortest path connecting~$u$
and $v$.

A \emph{tree} is a connected graph which does not contain any cycles.
We will often use the fact that for any subtree $T'$ of a tree $T$ and any
vertex $x \in T$ there is a unique vertex $y \in T'$ which minimises $d(x, y)$
over all $y \in T'$.
For any vertex $x \in T$ and edge $e \in E(T)$ incident to~$x$, the \emph{weight
of $e$ from $x$}, denoted $w_e(x)$, is the number of vertices $y \neq x$ of $T$
for which~$e$ is the first edge of the path from $x$ to $y$. Each vertex $y \neq
x$ of $T$ contributes to the weight from $x$ of precisely one edge incident to
$x$, so the sum of the weights from $x$ over all edges incident to $x$ is
$|T|-1$. Also, if $xy$ is an edge of $T$, then $w_e(x)+w_e(y) = |T|$.

A \emph{rooted tree} is a tree with a specified vertex $r$ as a \emph{root}. In
a rooted tree every vertex~$x$ other than the root has a \emph{parent}; this is
defined to be the unique neighbour $y$ of $x$ with $d(y, r) < d(x, r)$. If $y$
is the
parent of $x$ then we say that $x$ is a \emph{child} of $y$. A \emph{leaf} in a
tree is a vertex of degree one; so every vertex other than the root is a child
of some vertex, and every vertex apart from a leaf is a parent of some vertex.
An \emph{ancestral ordering} of the vertices of a tree is a linear order in
which the
root appears first and every other vertex appears after its parent.

A \emph{directed graph} $G = (V,E)$, or digraph, is formed by a vertex set $V$
and a set of edges~$E$, where every edge $e \in E$ is an ordered pair $(u, v)$ of
vertices of $G$. For $u, v \in V$ we write $u \rightarrow v$ or $v \leftarrow u
$ if $(u, v) \in
E(G)$. Also, for any vertex $v$ of $G$, $N^+(v)$ denotes the set of vertices $u$
such that $(v, u) \in E(G)$, and $N^-(v)$ denotes the set of vertices $u$ such
that $(u, v) \in E(G)$. $d^+(v)$ and $d^-(v)$ denote $|N^+(v)|$ and $|N^-(v)|$
respectively, and $\delta^+(G)$ and $\delta^-(G)$ are then defined to be the
minimum of $d^+(v)$
and $d^-(v)$ respectively over all vertices $v \in G$. The minimum
\emph{semidegree} is $\delta^0(G) = \min \{\delta^+(G), \delta^-(G)\}$. 
A \emph{tournament} on $n$
vertices is a digraph $G$ on~$n$ vertices in which for any distinct $u, v \in
V(G)$ precisely one of $u \rightarrow v$ and $u \leftarrow v$ holds. So a
tournament can be thought of as an orientation of the complete graph on~$n$
 vertices. Given a digraph $G$, the \emph{underlying graph} $G_{\textrm{under}}$ is the
graph on $V(G)$ in which there is an edge between~$u$ and $v$ if and only if
either $u \rightarrow v$ or $u \leftarrow v$. We define the distance $d(u, v)$
between distinct vertices $u$ and $v$ of a digraph $G$ to be the distance
between those two vertices in the underlying graph $G_{\textrm{under}}$. Also, if $G$ is a graph or digraph, and $H$ is a subgraph of $G$, then we write $G - H$ to denote $G[V(G) \sm V(H)]$, that is, the subgraph of $G$ induced on those vertices not in $H$.

A directed tree is a digraph $T$ for which the underlying graph $T_{\textrm{under}}$ is a
tree and in which at most one of $x \rightarrow y$ and $x \leftarrow y$ holds
for any pair of vertices $x$ and $y$ of $T$. We use the notation $x \rightarrow
y$ to distinguish a directed edge from an undirected edge, for which we use the
notation $xy$. Given a specified vertex $r$ as a
root, we define parents and children of vertices of the directed tree $T$
exactly
as in the underlying tree $T_{\textrm{under}}$. Similarly $\Delta(T) = \Delta(T_{\textrm{under}})$, and the
weight $w_e(x)$ of an edge
$e$ incident to a vertex $x$ is defined as in $T_{\textrm{under}}$. Also, for each vertex
$v$, $w^+(x)$ is the sum of $w_e(x)$ over all edges $e$ incident to $x$ directed
away from $x$, and $w^-(x)$ is the sum of $w_e(x)$ over all edges $e$ incident
to $x$ directed towards $x$. More generally, for a subtree $T'$ of $T$, $w^+(T')$ is
the sum of $w_e(x)$ over all edges $e$ directed from a vertex $x$ of $T'$ to a
vertex of $T - T'$, and $w^-(T')$ is the sum of $w_e(x)$ over all edges $e$
directed from a vertex of $T-T'$ to a vertex $x$ of $T'$. We say that a vertex of a digraph is a \emph{sink} vertex if it has no outneighbours, and a \emph{source} vertex if it has no inneighbours. Since a directed tree on $n$ vertices has $n-1$ edges, any directed tree must contain at least one sink vertex and at least one source vertex.

Throughout the paper we use the notation $x \ll y$ to
indicate that for any $y > 0$ there exists $x_0 >0$ such that for any $0 < x
\leq x_0$ the subsequent statements hold. Such statements with more variables
are defined similarly. Also, we will sometimes write `let $x \ll y$' when $y$
has an already fixed positive value; by this we mean that there exists some $x_0
> 0$ such that for any $0 < x < x_0$ the subsequent statements hold.
When we use asymptotics such as $o(1)$ we mean that these hold as $n\to \infty$
and all the other parameters are fixed.

\subsection{Probabilistic estimates}

The next lemma, relating to binomial distributions, will be used to show that in
the randomised algorithm we use in Section~\ref{sec:bounded_robust_exp},
the cluster to which a vertex is allocated is
almost independent of the cluster to which a vertex far away is allocated. We
use $\mathcal{B}(n, p)$ to denote the binomial distribution with parameters $n$
and $p$, that is, the number of successes in $n$ independent trials, each of
which has probability $p$ of success. So $\mathbb{E}(\mathcal{B}(n,p)) = np$.

\begin{lemma} \label{binom_uniform}
Suppose that $1/k \ll p, (1-p), \eps$, that $n \geq k^3/6$, and that $X =
\mathcal{B}(n,p)$. Then
for any $0 \leq r \leq k-1$, $$\mathbb{P}(X \equiv r \mod k) = (1
\pm \eps)/k.$$
\end{lemma}

\proof For each $x \in \{0, \dots, n\}$ let $p_x$ denote $\mathbb{P}(X = x)$, so
$p_x = \binom{n}{x}p^x(1-p)^{n-x}$. Let $\mu = np$, so $\mathbb{E}(X) = \mu$,
and let $p_\mu = \max\{p_{\lfloor \mu \rfloor}, p_{\lceil \mu \rceil}\}$, so
$p_x \leq p_\mu$ for any $x$. Moreover, if $x \leq y \leq \mu$ or $\mu \leq y
\leq x$ then $p_x \leq p_y$. So for any $r, i \in [k]$, \begin{align*}
\mathbb{P}(X \equiv r \mod k) &=
\sum_{\substack{0 \leq x \leq \mu - k \\ x\equiv r \mod k}} p_x +
\sum_{\substack{\mu - k < x \leq \mu + k \\ x\equiv r \mod k}} p_x+
\sum_{\substack{\mu + k < x \leq n \\ x\equiv r \mod k}} p_x
\\ &\leq \sum_{\substack{0 \leq x \leq \mu - k \\ x\equiv r \mod k}} p_{x+i} +
2 p_\mu +
\sum_{\substack{\mu + k < x \leq n \\ x\equiv r \mod k}} p_{x-k+i}
\\ &\leq \mathbb{P}(X \equiv r+i \mod k) +2 p_\mu.
\end{align*}

So $\mathbb{P}(X \equiv r \mod k) = 1/k \pm 2p_\mu = (1 \pm \eps)/k$ for any $r
\in [k]$, using a standard result (e.g. \cite{Bbook}, Section 1.2) on the
binomial distribution which states that $p_\mu = O(n^{-1/2}) = O(k^{-3/2})$.
\endproof

The following two results give useful tail estimates for random
variables. The first is an Azuma-type inequality which bounds the sum of many
small and almost independent random variables. This is derived in \cite{SV} from
a result in~\cite{Mc}. (\cite{SV} uses a random walk to embed trees in
sparse undirected graphs.) The second gives standard Chernoff-type bounds for the
binomial and hypergeometric distributions. The hypergeometric random variable
$X$ with parameters $(n,m,k)$ is defined as follows. Let $N$ be a set of size
$n$, and fix a set $S \subseteq N$ of size $|S|=m$. Now choose a set $T
\subseteq N$ of size $|T|=k$ uniformly at random. Then $X=|T \cap S|$. Note that
$\mathbb{E}X = km/n$.

\begin{lemma} [\cite{SV}, Proposition 1.1] \label{azuma}
Let $X_1, \dots, X_n$ be random variables taking values in $[0,1]$ such that for
each $k \in [n]$,
$$\mathbb{E}(X_k \mid X_1, \dots, X_{k-1}) \leq a_k.$$
Let $\mu \geq \sum_{i=1}^n a_i$. Then for any $0 < \delta < 1$,
$$\mathbb{P}(\sum_{i=1}^n X_i > (1+\delta) \mu) \leq e^{-\frac{\delta^2
\mu}{3}}. $$
\end{lemma}

\begin{prop} [\cite{JLR}, Corollary 2.3 and Theorem 2.10]\label{chernoff}
Suppose $X$ has binomial or hypergeometric distribution and $0<a<3/2$. Then
$\mathbb{P}(|X - \mathbb{E}X| \ge a\mathbb{E}X) \le 2
e^{-\frac{a^2}{3}\mathbb{E}X}$.
\end{prop}

\subsection{Regularity and Robust Outexpanders.}

To prove Theorem~\ref{main} we shall make use of a directed version of
Szemer\'edi's Regularity lemma. For this, we make the following definitions. If
$G$ is an undirected bipartite graph with vertex classes $X$ and $Y$, then the
density of $G$ is defined as
$$ d(X, Y) := \frac{e(X,Y)}{|X||Y|}.$$
Now, for any $\eps >0$, we say that $G$ is \emph{$\eps$-regular} if for any $X'
\subseteq X$ and $Y' \subseteq Y$ with $|X'| \geq \eps |X|$ and $|Y'| \geq \eps
|Y|$ we have $|d(X',Y') - d(X, Y)| < \eps$.

Given disjoint vertex sets $X$ and $Y$ in a digraph $G$, we use $G[X \rightarrow
Y]$ to denote the edges of $G$ directed from $X$ to $Y$. We say $G[X \rightarrow
Y]$ is \emph{$\eps$-regular with density $d$} if the underlying bipartite graph
of $G[X \rightarrow Y]$ is $\eps$-regular and has density $d$. Next we state
the degree form of the regularity lemma for digraphs. A regularity lemma for
digraphs was proven by Alon and Shapira~\cite{AS}. The degree form follows from
this in the same way as the undirected version (see~\cite{BCCsurvey}
for a sketch of the latter).

\begin{lemma}[Regularity Lemma for digraphs] \label{regularity_lemma}
For any $\eps, M'$ there exist $M, n_0$ such that if~$G$ is a digraph on $n \geq
n_0$ vertices and $d \in [0, 1]$, then there exists a partition of $V(G)$ into
$V_0, \dots, V_k$ and a spanning subgraph $G'$ of $G$ such that
\begin{itemize}
\item[(1)] $M' \leq k \leq M$,
\item[(2)] $|V_0| \leq \eps n$,
\item[(3)] $|V_1| = \dots = |V_k| =: m$,
\item[(4)] $d^+_{G'}(x) > d^+_{G}(x) - (d+\eps) n$ for all vertices $x \in V(G)$,
\item[(5)] $d^-_{G'}(x) > d^-_{G}(x) - (d+\eps) n$ for all vertices $x \in V(G)$,
\item[(6)] for all $i \in [k]$ the digraph $G'[V_i]$ is empty,
\item[(7)] for all $1 \leq i, j \leq k$ with $i \neq j$ the pair $G'[V_i
\rightarrow V_j]$ is $\eps$-regular and either has density 0 or density at least
$d$.
\end{itemize}
\end{lemma}

We refer to $V_1, \dots, V_k$ as \emph{clusters}.
Given a graph $G$ on $n$ vertices, we form the \emph{reduced digraph $R$ of $G$
with parameters $\eps, d$ and $M'$} by applying the regularity lemma with these
parameters to obtain $V_0, \dots, V_k$. $R$ is then the digraph on
vertex set $\{1, \dots, k\}$, with $i \rightarrow j$ an edge precisely when
$G'[V_i \rightarrow V_j]$ is $\eps$-regular with density at least $d$.

One particular regular structure will appear frequently in
Section~\ref{sec:bounded_robust_exp} and Section~\ref{sec:unbounded_robust_exp}.
We say that a digraph $G$ is an \emph{$\eps$-regular $d$-dense cycle of cluster
tournaments} if $V(G) = V_1 \cup \dots \cup V_k$, where the sets $V_i$ are
pairwise disjoint and of equal size\COMMENT{Check that this is always the
case!}, and for each $i$, $G[V_i]$ is a tournament and $G[V_i \rightarrow
V_{i+1}]$ is $\eps$-regular with density at least $d$ (where here and throughout
this paper addition and subtraction on the indices of clusters is to be taken
modulo $k$). We shall often refer to the sets $V_i$ as clusters, as we will
obtain
them by an application of the regularity lemma.

Now, let $V_1, \dots, V_k$ be disjoint sets of $m$ vertices, and let $G$ be a
digraph on vertex set $V_1 \cup \dots \cup V_k$. Let $S$ be a subset of some
cluster $V_i$. Then we say that $S$ is ($c, \gamma$)-\emph{good} if for any
$V_{i-1}'\subseteq V_{i-1}$ and $V_{i+1}' \subseteq V_{i+1}$ with $|V_{i-1}'|
\geq cm$ and $|V_{i+1}'| \geq cm$, $S$ contains at least $\gamma
\sqrt{m}$ vertices which each have at least $\gamma m$ inneighbours in
$V_{i-1}'$ and at least $\gamma m$ outneighbours in $V_{i+1}'$. Our main tool in
the use of regularity will be the next lemma, which states that if $G$ is a
regular and dense cycle of cluster tournaments, then any subset $V_i'$ of any
cluster $V_i$ with $|V_i'| \geq \gamma m/2$ contains a $(c, \gamma)$-good subset $S$
of size at most~$\sqrt{m}$.

\begin{lemma} \label{sqrt_sets}
Suppose that $1/m \ll \eps \ll \gamma \ll c, d$. Let $G$ be an $\eps$-regular
$d$-dense cycle of cluster tournaments on clusters $V_1, \dots, V_k$, each of
size $m$. Then for any $i$ and for any $V'_i \subseteq V_i$ of size $|V'_i| =
\gamma m/2$, there exists a $(c, \gamma)$-good set $S \subseteq V'_i$ with $|S|
\leq \sqrt{m}$.
\end{lemma}

\proof
Given $V'_i \subseteq V_i$ of size $|V'_i| = \gamma m/2$, choose $S \subseteq
V_i'$ at random by including
each vertex of $V_i'$ in $S$ with probability $1/\gamma\sqrt{m}$, independently
of the outcome for each other vertex. Then by Proposition~\ref{chernoff}, with
probability $1-o(1)$, $|S| \leq \sqrt{m}$.

Now, $G[V_{i-1} \rightarrow V_i']$ and $G[V_i' \rightarrow V_{i+1}]$ are each
$(2\eps/\gamma)$-regular with density at least $d/2$. So all but at most $2\eps
m/\gamma$ vertices $v_{i-1} \in V_{i-1}$ have at least $\gamma dm/5$
outneighbours in $V_i'$. Fix any such $v_{i-1} \in V_{i-1}$. Then  $G[V_i' \cap
N^+(v_{i-1}) \rightarrow V_{i+1}]$ is $(5\eps/\gamma d)$-regular with density at
least $d/2$. So all but at most $5\eps m/\gamma d$ vertices $v_{i+1} \in
V_{i+1}$ have at least $\gamma d^2m/20$ inneighbours in $V_i' \cap
N^+(v_{i-1})$. We therefore conclude that all but at most $7 \eps m^2/\gamma d$
pairs $(v_{i-1}, v_{i+1})$ with $v_{i-1} \in V_{i-1}$, $v_{i+1} \in V_{i+1}$
have at least $\gamma d^2m/20$ common neighbours in $V_i'$.

By Proposition~\ref{chernoff}, for each such pair
$(v_{i-1}, v_{i+1})$ the probability that $(v_{i-1}, v_{i+1})$ has fewer than
$d^2 \sqrt{m}/25$ common neighbours in $S$ decreases exponentially with~$m$,
whilst the number of such pairs is quadratic in $m$. Thus with probability
$1-o(1)$ our randomly selected $S$ will have the property that all but at most
$7 \eps m^2/\gamma d$ pairs $(v_{i-1}, v_{i+1})$ with $v_{i-1} \in V_{i-1}$,
$v_{i+1} \in V_{i+1}$ have at least $d^2 \sqrt{m}/25$ common neighbours in $S$.
We may therefore fix an outcome of our random choice of $S$ such that both of
these events of probability $1-o(1)$ occur.

So if $|V_{i-1}'| \geq cm$ and $|V_{i+1}'| \geq cm$, then we know that at least
$c^2m^2/2$ pairs $(v_{i-1}, v_{i+1})$ with $v_{i-1} \in V_{i-1}'$, $v_{i+1} \in
V_{i+1}'$ have at least $d^2 \sqrt{m}/25$ common neighbours $s \in S$. Thus
there are at least $c^2 d^2 m^{5/2}/50$ triples of such vertices $(v_{i-1}, s,
v_{i+1})$, so at least $c^2 d^2 \sqrt{m}/100 \geq \gamma \sqrt{m}$ vertices in
$S$ must lie in the common neighbourhood of at least $c^2 d^2 m^2/100$ such
pairs $(v_{i-1}, v_{i+1})$. (Otherwise there would be fewer than $|V'_{i-1}||V'_{i+1}|c^2d^2\sqrt{m}/100 + |S|c^2d^2m^2/100
\leq c^2d^2m^{5/2}/50$ such triples $(v_{i-1},s,v_{i+1})$.)  Each of these vertices therefore has at least $c^2
d^2 m/100 \geq \gamma m$ neighbours in each of $V_{i-1}'$ and $V_{i+1}'$, as
required.\endproof

We will also make use of the following well known observation, which says that
if $G$ is a
regular and dense cycle of cluster tournaments on clusters $V_1, \dots, V_k$,
and we select subsets $U_1 \subseteq V_1, \dots, U_k \subseteq V_k$ uniformly at
random, then with high probability the restriction of $G$ to these subsets is
also regular and dense. This follows from a lemma of Alon et al.~\cite{Aetal}
showing that
$\eps$-regularity is equivalent to almost all vertices having the expected degree
and almost all pairs of vertices having the expected common neighbourhood size. We
include the proof for completeness.

\begin{lemma} \label{regular_restriction}
Suppose that $1/m \ll k \ll \eps \ll \eps' \ll d$ and that $m^{1/3} \leq m' \leq
m$. Let $G$ be an $\eps$-regular $d$-dense cycle of cluster tournaments on
clusters $V_1, \dots, V_k$, each of size $m$. For each $i \in [k]$, choose $U_i
\subseteq V_i$ of size $m'$ uniformly at random, and independently of all other
choices. Then with probability $1-o(1)$, $G[U_1 \cup \dots \cup U_k]$ is an
$\eps'$-regular $d/2$-dense cycle of cluster tournaments.
\end{lemma}

\proof
We need to show that with high probability, $G[U_i \rightarrow U_{i+1}]$ is
$\eps'$-regular with density at least $d/2$ for each $i$. So fix some $i \in
[k]$, and let $d_i \geq d$ be the density of $G[V_i \rightarrow V_{i+1}]$. Also,
let $B_i$ be the set of vertices $v \in V_i$ for which $|N^+(v) \cap V_{i+1}|
\neq (d_i \pm \eps)m$, and let $D_i$ be the set of pairs $v_1 \neq v_2$ of
vertices of $V_i$ for which $|N^+(v_1) \cap N^+(v_2) \cap V_{i+1}| \neq (d_i^2 \pm
3\eps)m$. Then since $G[V_i \rightarrow V_{i+1}]$ is $\eps$-regular, $|B_i|
\leq 2\eps m$. Also, there are at most $2\eps m^2$ pairs in $D_i$ which contain
a vertex of $B_i$, and each $v \in V_i \sm B_i$ lies in at most $2\eps m$ pairs
in $D_i$, so $|D_i| \leq 4 \eps m^2$. So let $B_i' = B_i \cap U_i$ and similarly
let $D_i'$ consist of the pairs in $D_i$ for which both vertices of the pair are
in $U_i$. Then by Proposition~\ref{chernoff}, the probability that either
$|B_i'| > 4\eps m'$ or $|D_i'| > 8\eps (m')^2$ declines exponentially with $m$.

Now, for each of the at most $m'$ vertices $v \in U_i \sm B_i$, by
Proposition~\ref{chernoff} the probability that $|N^+(v) \cap U_{i+1}| \neq (d_i
\pm 2\eps)m'$ decreases exponentially with $m$. Also, for each of the at most
$\binom{m'}{2}$ pairs $v_1 \neq v_2$ with $v_1, v_2 \in U_i\setminus D_i$,
the probability that $|N^+(v_1) \cap N^+(v_2) \cap U_{i+1}| \neq (d_i^2 \pm
4\eps) m'$ decreases exponentially with $m$. So with probability $1-o(1)$, for
each $i$ none of these events of exponentially declining probability will
hold.

Fix such an outcome of our random choices. Then for each $i$ at least $(1-4\eps)
m'$ vertices $v \in U_i$ have $|N^+(v_1) \cap U_{i+1}| = (d_i \pm 2\eps)m'$ and
at least $\binom{m'}{2} - 8 \eps (m')^2$ pairs $v_1, v_2 \in U_i$ have
$|N^+(v_1) \cap N^+(v_2) \cap U_{i+1}| = (d_i^2 \pm 4\eps)m'$. It then
immediately follows from Lemma~3.2 of~\cite{Aetal} that for each $i$, $G[U_i
\rightarrow U_{i+1}]$ is $\eps'$-regular (and it is clear that this has density
at least $d_i/2 \geq d/2$), as desired.
\endproof

We now turn to the concept of a robust outexpander. Let $\mu > 0$, let $G$ be a
digraph on~$n$ vertices, and let $S \subseteq V(G)$. Then the \emph{robust
$\mu$-outneighbourhood} of $S$, denoted $RN^+_\mu(S)$, is defined to be the set
of vertices of $G$ with at least $\mu n$ inneighbours in $S$. For constants $0 <
\mu \leq \nu < 1$, we say that a digraph $G$ on $n$ vertices is a \emph{robust
$(\mu, \nu)$-outexpander} if $|RN^+_\mu(S)| \geq |S| + \mu n$ for all $S
\subseteq V(G)$ with $\nu n< |S| < (1-\nu)n$. A recent result from~\cite{KOT}
(which in turn relies on results from~\cite{KKO, KeeKO}) states that every robust
outexpander with linear minimum semidegree contains a Hamilton cycle.
We shall make use of this to prove the next lemma, which states that if a
tournament $G$ is a robust outexpander then $G$ contains a regular and dense
cycle of cluster tournaments which covers almost all of the vertices of $G$.
We will use this structure when we embed a tree $T$ in a tournament $G$
which is a robust outexpander.

\begin{lemma} \label{KOT_combined}
Suppose that $1/n \ll 1/M \ll 1/M' \ll \eps \ll d \ll \mu \ll \nu \ll \eta<1$.
Let~$G$ be a tournament on $n$ vertices which is a robust $(\mu,
\nu)$-outexpander with $\delta^0(G) \geq \eta n$. Then~$G$ contains an
$\eps$-regular $d$-dense cycle of cluster tournaments on clusters $V_1, \dots,
V_k$, where $|\bigcup_{i=1}^k V_i| > (1-\eps)n$, and $M' \leq k \leq M$.
\end{lemma}

\proof
Let $R$ be the reduced digraph of $G$ with parameters $\eps, d$ and $M'$ obtained
by applying Lemma~\ref{regularity_lemma}, and let $k = |R|$, so $M' \leq k \leq
M$. Then by Lemma~12 of \cite{KOT}, $\delta^0(R) \geq \eta|R|/2 $, and $R$ is a
robust $(\mu/2, 2\nu)$-outexpander. Then by Theorem~14 of \cite{KOT}, which
states that any robust outexpander of linear minimum semidegree contains a Hamilton
cycle, we know that $R$ contains a Hamilton cycle. Let $V_1, \dots, V_k$ be the
clusters of $R$ in the order of the cycle. Then $|\bigcup_{i=1}^k V_i| >
(1-\eps)n$, $G[V_i]$ is a tournament for each $i$ and since $V_i \rightarrow
V_{i+1}$ is an edge of $R$ for each $i$, $G'[V_i \rightarrow V_{i+1}]$ is
$\eps$-regular with density at least $d$. (Here $G'$ is the spanning subgraph of $G$
obtained by Lemma~\ref{regularity_lemma}.)
\endproof

Of course, we will sometimes need to embed a tree $T$ in a tournament $G$ which
is not a robust outexpander. The next lemma will be a useful tool in this
situation; it states that if a tournament $G$ is not a robust outexpander then
$V(G)$ can be partitioned into two sets so that most edges between the two sets
have the same direction.

\begin{lemma} \label{not_robust_expander_split}
Suppose that $1/n \ll \mu \ll \nu$, that $G$ is a tournament on $n$ vertices
and that $G$ is not a robust $(\mu, \nu)$-outexpander. Then we can partition
$V(G)$ into sets $S$ and $S'$ such that $\nu n < |S|, |S'| < (1-\nu)n$
and $e(G[S \rightarrow S']) \leq 4 \mu n^2$.
\end{lemma}

\proof
Since $G$ is not a robust $(\mu, \nu)$-outexpander there exists $S \subseteq
V(G)$
such that $|RN_\mu^+(S)| < |S| + \mu n$ and $\nu n < |S| < (1-\nu)n$.
Choose such an $S$, and let $S' = V(G) \sm S$, so $\nu n < |S'| <
(1-\nu)n$ also.

Since $G$ is a tournament, at most $2 \mu n+1$ vertices $v \in S$ have
$d^-_{G[S]}(v) < \mu n$, and so at most $2 \mu n+1$ vertices $v \in S$ have $v
\notin RN^+_\mu(S)$. So $|RN_\mu^+(S) \sm S| \leq 3 \mu n+1$, and so the number of
edges from $S$ to $S'$ is at most\COMMENT{Since there are at most $(3 \mu n+1) |S|$
edges from $S$ to $|RN_\mu^+(S) \sm S|$, and every vertex of $S'\sm RN_\mu^+(S)$
has at most $\mu n$ inneighbours in $|S|$, adding at most another $\mu n |S'|$
edges from $S$ to $S'$.} $(3 \mu n+1) |S| + \mu n |S'| \leq 4 \mu n^2$.
\endproof

\subsection{Basic Tree Properties}

In this section, we shall prove several lemmas which we shall make use of in
proving Theorem~\ref{main}. The first two of these will enable us to split a
tree into several pieces with properties that will be useful for the analysis
of the randomised embedding algorithm used in Section~\ref{sec:bounded_robust_exp}.

\begin{lemma} \label{single_tree_split}
Let $T$ be a tree on $n \geq 3$ vertices. Then there exist subtrees
$T'$ and $T''$ of $T$ such that $T'$ and $T''$ intersect in precisely one vertex
of $T$, every edge of $T$ lies in precisely one of $T'$ and $T''$, and $e(T'),
e(T'') \geq e(T)/3$.
\end{lemma}

\proof
We begin by showing that $T$ must contain a vertex $v$ such that every edge $e$
incident to $v$ has $w_e(v) \leq n/2$. Recall that if $e = uv$, then
$w_e(u)+w_e(v) = n$, and so at most one of $w_e(u) > n/2$ and $w_e(v) > n/2$ can
hold. Since $T$ contains $n$ vertices and $n-1$ edges, by the pigeonhole
principle $T$ contains a vertex $v$ so that no edge $e$ incident to $v$ has
$w_e(v) > n/2$.

Now, choose such a vertex $v$ in $T$, and let $v_1, \dots, v_r$ be the
neighbours of $v$ in $T$. For each $i$, let $S_i$ be the set of vertices $x$ of
$T$ such that $v_i$ lies on the path from $v$ to $x$. Then every vertex of $T$
other than $v$ lies in precisely one set $S_i$. Now, for each $i$, let $T_i$ be
the tree $T[S_i \cup \{v\}]$. Then each $T_i$ is a subtree\COMMENT{
Clearly $T_i$ is acyclic, and there are no edges between $S_i$ and $S_j$ for $i
\neq j$. So $T[S_i]$ must be connected as $T$ is connected, and $vv_i$ is an
edge so $T[S_i \cup \{v\}]$ is connected.
} of $T$ and every edge of $T$ is contained in precisely one\COMMENT{
If $v \in e$ then $e$ is of the form $vv_i$ and lies in $T_i$ and no other
$T_j$; if $v \notin e$ then $u, v$ lie in the same set $S_i$ and hence lie in
this $T_i$ and no other.
} $T_i$. So $\{e(T_i): i \in [r]\}$ is a set of positive integers, none
greater\COMMENT{
Since if $e(T_i) > 2(n-1)/3$ then $e(T_i) > n/2$ so $w_{vv_i}(v) > n/2$,
contradiction.
} than $2(n-1)/3$, which sum to $n-1$. Thus there exists\COMMENT{
If $a_1, \dots, a_r$ are non-negative integers with $a_1+\dots+a_r = m$ and $a_i
\leq 2m/3$ for each $a_i$, then for some $A \subseteq [r]$ we have $m/3 \leq
\sum_{i \in A} a_i \leq 2m/3$. To see this, note that if some $a_i \geq m/3$,
then we are done (take $A = \{i\}$). If not, then let $A$ be a minimal (w.r.t.
inclusion) subset of $[r]$ with $\sum_{i \in A} a_i \geq m/3$.
} $A \subseteq [r]$ such that the sum of elements of $\{e(T_i): i \in A\}$ lies
between $(n-1)/3$ and $2(n-1)/3$. Then if we take $T' = \bigcup_{i \in A} T_i$
and $T'' = \bigcup_{i \notin A} T_i$ then $T'$ and $T''$ satisfy the conditions
of the lemma (in particular, $T' \cap T'' = \{v\}$).
\endproof

\begin{lemma} \label{split_trees}
Suppose that $1/n \ll 1/\Delta, \eps, 1/k$. Let $T$ be a tree on $n$ vertices
satisfying $\Delta(T) \leq \Delta$ and rooted at $t_1$. Then there
exist pairwise disjoint subsets $F_1, \dots, F_r$ of $V(T)$, and vertices $v_1,
\dots, v_r$ (not necessarily distinct) of $T$ such that:
\begin{itemize}
\item[(1)] $|\bigcup_{i \in [r]} F_i| \geq (1-\eps)n $.
\item[(2)] $|F_i| \leq n^{2/3}$ for each $i$.
\item[(3)] For any $i \in [r]$, let $x \in \{t_1\} \cup \bigcup_{j<i} F_j$, and
let $y \in F_i$. Then the path from $x$ to $y$ in~$T$ includes the vertex $v_i$.
\item[(4)] For any $y \in F_i$ we have $d_T(v_i, y) \geq k^3$.
\end{itemize}
\end{lemma}

\proof We begin by splitting $T$ into a family $\F$ of subtrees of $T$ by
repeated use of Lemma~\ref{single_tree_split}. So initially let $\F = \{T\}$,
and then we repeat the following step. Let $T_{large}$ be the largest member of
$\F$. Use Lemma~\ref{single_tree_split} to split $T_{large}$ into
subtrees $T'$ and $T''$ which intersect in a single vertex, partition the edges
of $T_{large}$, and satisfy $e(T'), e(T'') \geq e(T_{large})/3$. Then remove
$T_{large}$ from $\F$, and replace it by the two smaller trees $T'$ and $T''$.
After at most $3n^{1/3}$ steps\COMMENT{
After $3n^{1/3}$ steps we have $|\F| = 3n^{1/3}$ and so some $T_i \in \F$ has
$|T_i| \leq n^{2/3}/3$. Then $|T_i|$ must have been formed by splitting a
subtree $T'$ with $|T'| \leq n^{2/3}$, which must have been the largest tree in
$\F$ at that point, so at that point, every $T \in \F$ had $|T| \leq n^{2/3}$.
} we must have that $|T^*| \leq n^{2/3}$ for every $T^* \in \F$. At this point
we terminate the process.

Observe that if $T', T'' \in \F$, then $T'$ and $T''$ intersect in at most one
vertex.\COMMENT{By induction; this is clearly true when $\F = \{T\}$, and when
some $T \in \F$ is split into $T'$ and $T''$ then since every $T_i \in \F$
intersected $T$ in at most one vertex they intersect $T'$ and $T''$ in at most
one vertex.}
Now, form a graph $\G_\F$ with vertex set $\F$ and with an edge between $T', T''
\in \F$ if and only if $T'$ and $T''$ have a common vertex. Then $\G_\F$ is
connected\COMMENT{
If not, then we can write $\F = \F' \cup \F''$ such that $T' \cap T'' =
\emptyset$ for any $T' \in \F'$ and $T'' \in \F''$. Since every edge of $T$ lies in
a member of $\F$, this yields non-empty sets of vertices of $T$ with no edge
between them, contradicting $T$ being connected.
}, and so contains a spanning tree $\T_\F$. Choose~$T_0$ to be a member of $\F$
containing the root $t_1$ of $T$, and let $T_0, T_1, \dots, T_r$ be an ancestral
ordering of the members of $\F$ (thought of as vertices of the tree $\T_\F$).
Now, for each $1 \leq i \leq r$ let $v_i$ be the common vertex of $T_i$ and its
parent in $\T_\F$. Then define $F_i$ for each $i \in [r]$ by
$$F_i = V(T_{i}) \sm \{x \in T: d_T(v_i, x) < k^3\}.$$
It remains to show that $F_1, \dots, F_r$ and $v_1, \dots, v_r$ satisfy the
required properties. (4) is immediate from the definition of $F_i$, and (2)
holds since each $T_i$
contained at most $n^{2/3}$ vertices. For (1), observe that every vertex of $T$
was contained in at least one of the subtrees $T_i$, and that in forming the
sets $F_i$, we deleted at most $\Delta^{k^3}$ vertices from each of the at most
$3n^{1/3}$ sets $V(T_i)$, so in total at most $3n^{1/3} \Delta^{k^3} \leq \eps
n$ vertices of $T$ are not contained in any of the sets $F_i$.

For condition (3), suppose that $T_1'T_2'T_3'$ is a path in $\T_\F$, and some
vertex $v$ lies in $T_1' \cap T_3'$, but $v \notin T_2'$. Let $v' \in T_1' \cap
T_2'$ and let $v'' \in T_2' \cap T_3'$. Then $v' \neq v''$, as otherwise $T_1'$
and $T_3'$ would have a common vertex other than $v$. So there is a path from
$v'$ to $v''$ in $T$ which does not contain $v$, so $T$ contains a cycle, giving
a contradiction. Similarly it follows that for any path $T_{i_1}\dots T_{i_j}$
in $\T_\F$, if $T_{i_1}$ and $T_{i_j}$ have a common vertex $v$, then $v$ lies
in each of $T_{i_1}, \dots, T_{i_j}$, and so if $T_{i_{j-1}}$ is the parent of
$T_{i_j}$ in $\T_\F$ then $v = v_{i_j}$. Now, for any $i \in [r]$, if $x \in
\{t_1\} \cup \bigcup_{j<i} F_j$ and $y \in F_i$, then $x \in T_j$ for some $0
\leq j<i$ and $y \in T_i$. Let $T_jT_1'\dots T_s'T_i$ be the path from~$T_j$ to
$T_i$ in $\T_\F$, then $T_s'$ is the parent of $T_i$ in $\T_\F$.
So $T_j \cup T_1' \cup \dots \cup T_s'$ contains a path $P_1$ from~$x$ to $v_i$
and $T_i$ contains a path $P_2$ from $v_i$ to $y$. But the property we have
proved before implies that $P_1$ and $P_2$ only intersect in $v_i$. Thus $P_1
\cup P_2$ is the path in $T$ from $x$ to $y$, and $v_i \in P_1 \cup P_2$, as
required. It also follows that the sets $F_i$ are pairwise disjoint.\COMMENT{Since if
$F_j$ and $F_i$ with $j < i$ have a common vertex $v$, then $v \neq v_i$ (since
$v_i \notin F_i$), and so $v_i$ cannot be a common vertex of $T_j$ and $T_i$,
and so $v_i \notin T_j$. But then the path in $T_j$ from any vertex other than
$v$ to $v$ is a path in $T$ from a vertex of $T_j$ to a vertex of $T_i$, not
containing $v_i$, contradiction.}
\endproof

Recall that in Section~\ref{sec:bounded_robust_exp}, we will describe a randomised algorithm for
embedding the vertices of a tree $T$ in a digraph $G$. Whenever this algorithm
embeds a vertex $t$ of $T$ in $G$, it will reserve a set of vertices of $G$ in
which to embed the children of $t$. No other vertices may be embedded in this
set until all the children of $t$ have been embedded. For this to work, we need
to ensure that there will not be too many of these reserved sets at any point.
This motivates the following definition. If $T$ is a rooted tree on $n$
vertices, then we say that an ancestral ordering of the vertices of $T$ is
\emph{tidy} if it has the property that for any initial segment $\mathcal{I}$ of
the order, at most $\log_2 n$ vertices in $\mathcal{I}$ have a child not in
$\mathcal{I}$. The following lemma shows that such an order exists for any
tree~$T$.

\begin{lemma} \label{open_vertices}
Let $T$ be a tree on $n$ vertices rooted at some $t_0 \in T$. Then there exists
a tidy ancestral ordering of the vertices of $T$.
\end{lemma}

\proof We shall prove that for any $r$, the vertices of any rooted tree $T$ on
fewer than~$2^r$ vertices can be given an ancestral ordering so that fewer than
$r$ vertices in any initial segment~$\mathcal{I}$ have neighbours outside
$\mathcal{I}$. Indeed, suppose that this statement is false, and let $T$ rooted
at~$t_0$ be a counterexample of minimal order, say of order $n$. Let $r$ be
minimal such that $n <2^r$. Then let $T_1, \dots, T_s$ be the components of $T -
t_0$, ordered in increasing size, and let $t_i$ be the neighbour of $t_0$ in
$T_i$. We shall think of $t_i$ as the root of the tree $T_i$. Then $|T_i| <
2^{r-1}$ for $i \leq s-1$, and $T_s < 2^r$. So since $T$ was a minimal
counterexample, we can find an ancestral ordering of the vertices of each $T_i$
so that for any $i \leq s-1$, any initial segment of the order of the vertices
of $T_i$ contains fewer than $r-1$ vertices with children outside the initial
segment, and any initial segment of the order of the vertices of $T_s$ contains
fewer than $r$ vertices with children outside the initial segment. Now, we order
the vertices of $T$ as follows. Begin with~$t_0$, then add the vertices of $T_1$
in their order. Next, add the vertices of $T_2$ in their order, and continue in
this fashion. Since the order of the vertices of each $T_i$ was ancestral, this
order is also ancestral. Also, any initial segment $\mathcal{I}$ of this order
contains fewer than $r$ vertices with children outside $\mathcal{I}$,
contradicting the choice of $T$, and therefore proving the lemma.
\endproof

\section{Embedding trees of bounded maximum degree in a robust
outexpander.}\label{sec:bounded_robust_exp}

\subsection{Introduction}

Our aim in this section is the following lemma on embedding trees of bounded
maximum degree in robust outexpander tournaments.

\begin{lemma} \label{robust_expander_case_bounded_deg}
Suppose that $1/n \ll \mu \ll \nu \ll \eta \ll \alpha, 1/\Delta$, that
 $G$ is a
tournament on $(1+\alpha)n$ vertices which is a robust $(\mu,
\nu)$-outexpander
with $\delta^0(G) \geq \eta n$ and that $T$ is a directed tree on~$n$ vertices
with $\Delta(T) \leq \Delta$. Then $G$ contains a copy of $T$.
\end{lemma}

The proof of this lemma shows that we could actually put $1/\Delta$ lower down in
the hierarchy, but this is how will shall apply this lemma later on.
To prove this lemma, we begin by applying Lemma~\ref{KOT_combined} to find  a
regular and dense cycle of cluster tournaments in $G$, containing almost all of
the vertices of $G$. We will then use Lemma~\ref{cluster_cycle_bounded_deg} to
find a copy of $T$ within this structure. This
lemma is stated separately, and in a stronger form than necessary, as we shall
also make use of it in Section~\ref{sec:unbounded_robust_exp}.

\begin{lemma} \label{cluster_cycle_bounded_deg}
Suppose that $1/n \ll 1/k, 1/\Delta \ll \eps \ll d \ll \alpha \leq 2$, and
that $m =
n/k$. Let $G$ be an $\eps$-regular $d$-dense cycle of cluster tournaments on
clusters $V_1, \dots, V_k$ of equal size $(1+\alpha) m$. Let $v^*$ be a
vertex of $V_1$ with at least $d^2m$ inneighbours in $V_k$ and at least $d^2m$
outneighbours in $V_2$. Finally, let $T$ be a directed tree on $n$ vertices,
rooted at $t_1$ and with $\Delta(T) \leq \Delta$. Then~$G$ contains a
copy of $T$, where the vertex $t_1$ of $T$ corresponds to the vertex~$v^*$ of
$G$.
\end{lemma}

The main problem in achieving this is to allocate the vertices of $T$
to the clusters $V_i$ in such a way that we can then use the $\eps$-regularity
of each $G[V_i \rightarrow V_{i+1}]$ to embed the
vertices of $T$ in $G$. When we say we allocate
$v$ to $V_i$ this means that $v$ will be embedded to a vertex of $V_i$, but
this embedding has not been fixed yet. We wish to allocate each vertex of $T$
to a cluster $V_i$ so that for most edges $u \rightarrow v$ of $T$, if $u$ is
allocated to $V_i$ then $v$ is allocated to~$V_{i+1}$. So if $u$ is allocated
to a cluster $V_i$ and $u \rightarrow v$ then we say that the \emph{canonical
allocation} of~$v$ is to the cluster
$V_{i+1}$, whereas if $u \leftarrow v$ then we say that the canonical allocation
of~$v$ is to the cluster $V_{i-1}$. If we allocate $v$ canonically, then we say
that the edge between $u$ and~$v$ has been allocated canonically. One way of
allocating the vertices of $T$ to the clusters $V_i$ would be to begin by
allocating the root $t_1$ to $V_1$,
and then to allocate all remaining vertices canonically.
However, to successfully embed the vertices of $T$ within the clusters to which
they are allocated we
will need the vertices of $T$ to be approximately evenly distributed amongst the
$k$ clusters. This method will usually not achieve this, for example if $T$ is
an anti-directed path.

To obtain  an `even distribution' for any sufficiently large tree of bounded
maximum degree, we modify the method so that some vertices (selected
randomly) are allocated to the same cluster as their parent, rather than being
allocated canonically. However, having large components of vertices which are
allocated to the same cluster may prevent a successful embedding of these
vertices within this cluster, and so we shall also require that such components
are small. This is the motivation behind the Vertex Allocation Algorithm given
in the next subsection, which we shall use to allocate the
vertices of $T$.

\subsection{Allocating the vertices of $T$.}

We shall use the following random process to allocate the vertices of $T$ to
the clusters $V_i$.

\medskip
\noindent {\bf Vertex Allocation Algorithm:}

\emph{Input:} A directed tree $T$ on $n$ vertices, a root vertex $t_1 \in T$,
and clusters
$V_1, \dots, V_k$.

\emph{Initialisation:} Choose an ancestral ordering $t_1, \dots, t_n$ of the
vertices of $T$.

\emph{Procedure:} At time $\tau = 1$, allocate $t_1$ to $V_1$.
At time $\tau \geq 1$, we shall allocate $t_\tau$. Let $t_\sigma$ be the parent
of $t_\tau$, which must have appeared before $t_\tau$ in the ordering and has
therefore already been allocated. Then:
\begin{itemize}
\item If $d_T(t_\tau, t_1)$ is odd, then allocate $t_\tau$ canonically.
\item If $d_T(t_\tau, t_1)$ is even, then allocate $t_\tau$ to
the same cluster as $t_\sigma$ with probability $1/2$, and allocate $t_\tau$
canonically with probability $1/2$ (where these choices are made independently
for each vertex).
\end{itemize}

\emph{Termination:} Terminate when every vertex of $T$ has been processed and
therefore allocated to some cluster $V_j$.

\medskip

Note that the cluster to which a vertex $t$ is allocated by this algorithm
depends only on the
cluster to which its parent vertex was allocated and the outcome of the random
choice when
embedding $t$ (if $d(t, t_1)$ is even). Since these choices were independent, the
probability of any possible outcome does not depend on which ancestral order of
the vertices was chosen in the initialisation step. Now, we say that an edge of
$T$ is \emph{allocated within a cluster} if both of its endvertices are
allocated to the same cluster. Then we say that an allocation of the vertices of
a directed tree $T$ to clusters $V_1, \dots, V_k$ is \emph{semi-canonical} if
\begin{itemize}
\item[(i)] every edge of $T$ is either allocated canonically or is allocated
within a cluster,
\item[(ii)] every edge of $T$ incident to $t_1$ is allocated canonically, and
\item[(iii)] every component of the subgraph of $T$ formed by all edges
allocated within a cluster contains at most $\Delta(T)$ vertices.
\end{itemize}

The next lemma shows that if we allocate the vertices of a directed tree $T$ to
clusters $V_1, \dots, V_k$ by applying the Vertex Allocation Algorithm, then the
allocation obtained will be semi-canonical, and also that if vertices $t$ and
$t'$ are far apart in $T$ then the cluster to which $t$ is allocated is almost
independent of the cluster to which $t'$ is allocated. As a consequence, if~$T$
is sufficiently large and has bounded maximum degree, each cluster will have
approximately equally many vertices of $T$ allocated to it. These properties
will allow us to embed the vertices of such a $T$ into a regular and dense cycle
of cluster tournaments $G$ in the next subsection.

\begin{lemma} \label{allocation_props}
Let $T$ be a directed tree on $n$ vertices rooted at $t_1$. Allocate the
vertices of $T$ to clusters $V_1, \dots, V_k$ by applying the Vertex Allocation
Algorithm. Then the following properties hold.
\begin{itemize}
\item[(a)] The allocation obtained will be semi-canonical.
\item[(b)] Suppose that $1/k \ll \delta$. Let $u$ and $v$ be
vertices of $T$ such that $u$ lies on the path from $t_1$ to $v$, and $d_T(u, v)
\geq k^3$. Then for any $i, j \in [k]$, $$\mathbb{P}(v \textrm{ is allocated to
} V_i \mid u \textrm{ is allocated to } V_j) = \frac{1 \pm \delta/4}{k}.$$
\item[(c)] Now suppose also that $1/n \ll 1/\Delta, 1/k \ll \delta$, and that
$\Delta(T) \leq \Delta$. Then with probability $1-o(1)$, each of the $k$
clusters $V_i$ has at most $(1 + \delta)m$ vertices of $T$ allocated to it,
where $m = n/k$.
\end{itemize}
\end{lemma}

\proof
(a) The Vertex Allocation Algorithm allocates every vertex either canonically or
to the same cluster as its parent, so every edge will be allocated canonically
or within a cluster. Furthermore, a vertex $t$ can only be allocated to the same
cluster as its parent if $d(t_1, t)$ is even, and so each edge incident to $t_1$
is allocated canonically. Finally, since edges allocated within a cluster can
only be formed when we allocate $t_i$ such that $d(t_i, t_1)$ is even, any such
component is a star formed by some $t_j$ and some of the children
of~$t_j$.\COMMENT{Since $t_j$ has at most $\Delta(T)-1$ children, the result
follows.}

(b) Since the order in which the vertices are allocated is ancestral, at the
stage in our algorithm when we have just allocated $u$, no other vertex on the
path $P(u, v)$ in $T$ from $u$ to~$v$ has yet been allocated. So suppose that we
have just allocated $u$ to cluster $V_j$, let $\ell$ be the length of $P(u, v)$,
so $\ell \geq
k^3$, and let $u=v_0, v_1, \dots, v_{\ell}=v$ be the vertices of $P(u, v)$. Then
let $E = \{i \geq 1\colon d(v_i, t_1) \textrm{ is even}\}$, so $E$ indicates the
vertices
with a random element in their allocation, and let $O = [\ell] \sm E$, so $O$
indicates the vertices which are allocated deterministically. We then split the
edges of $P(u, v)$ into four classes:
\begin{align*}
F_{\textrm{canon}} &= \{v_{i-1} \rightarrow v_i: i \in O\}
\\B_{\textrm{canon}} &= \{v_{i-1} \leftarrow v_i: i \in O\}
\\F_{\textrm{random}} &= \{v_{i-1} \rightarrow v_i: i \in E\}
\\B_{\textrm{random}} &= \{v_{i-1} \leftarrow v_i: i \in E\}.
\end{align*}

Then every edge of $P(u, v)$ lies in one of these 4 sets, and so
$|F_{\textrm{canon}}|+|B_{\textrm{canon}}|+|F_{\textrm{random}}|+|B_{\textrm{random}}| = \ell$. Furthermore, each
edge in $F_{\textrm{canon}}$ will be allocated canonically, and hence from some $V_i$ to
$V_{i+1}$. Similarly, edges in $B_{\textrm{canon}}$ will be allocated from some $V_i$ to
$V_{i-1}$. Meanwhile, edges in $F_{\textrm{random}}$ or $B_{\textrm{random}}$ will be allocated
from some $V_i$ to
$V_{i+1}$ or $V_{i-1}$ respectively with probability $1/2$, and within some
$V_i$ with probability $1/2$. So let $R$ be the sum of the number of edges from
$F_{\textrm{random}}$ which are allocated canonically and the number of edges from
$B_{\textrm{random}}$ which are \emph{not} allocated canonically. Since the outcome of the
random experiment for each edge is independent of the outcome for any other edge, $R$
has distribution $\mathcal{B}(|E|, 1/2)$. Now,~$u$ was allocated to cluster
$V_j$, and so $v$ will be allocated to cluster
$V_i$, where
$$i \equiv j + |F_{\textrm{canon}}| - |B_{\textrm{canon}}| + R - |B_{\textrm{random}}| \mod k.
$$
But since $|E| \geq \lfloor \ell/2 \rfloor \geq k^3/3$, Lemma~\ref{binom_uniform} applied with $X=R$ and $n=|E|$ implies that for any $r \in [k]$, the probability that
$i = r$ is $\frac{1 \pm \delta/4}{k}$, as desired.

(c) Use Lemma~\ref{split_trees} to choose pairwise disjoint subsets $F_1, F_2,
\dots, F_r$ of $V(T)$ and vertices $v_1, \dots, v_r \in V(T)$ such that
$|\bigcup_{i \in [r]} F_i| \geq
(1-\delta/2k)n$ and $|F_i| \leq n^{2/3}$ for each $i$, also such that if $j <
i$, then any path
from $t_1$ or any vertex of $F_j$ to any vertex of $F_i$ passes through the
vertex $v_i$,
and finally such that $d(v_i, F_i) \geq k^3$. We shall prove that $(\dagger)$
with probability $1-o(1)$,
the total number of vertices from any of the sets $F_i$
allocated to cluster $V_j$ is at most $(1+\delta/2)m$. This will prove the
lemma, as the number of
vertices of $T$ not contained in any of the sets $F_i$ is at most $\delta m/2$,
and so in total at most $(1+\delta)m$ vertices of $T$ are allocated to any
cluster
$V_j$.

To prove $(\dagger)$, define random variables $X_i^j$ for each $i \in [r]$, $j
\in [k]$ by
$$X_i^j = \frac{\textrm{\# of vertices of } F_i \textrm{ allocated to cluster }
V_j}{n^{2/3}},$$
so that each $X_i^j$ lies in the range $[0,1]$.
Then since the cluster to which a vertex $t$ of $T$ is allocated is dependent
only on the cluster to which the parent of $t$ is allocated and on the outcome
of the random choice made when allocating $t$, $\mathbb{E}(X_i^j \mid X_{i-1}^j,
\dots, X_1^j, v_i \in V_s) = \mathbb{E}(X_i^j\mid~v_i\in V_s)$ for all $s \in
[k]$. Here we write $v_i \in V_s$ to denote the event that $v_i$ is allocated
to $V_s$. So for any $i$ and $j$,
\begin{align*}
\mathbb{E}(X_i^j \mid X_{i-1}^j, \dots, X_{1}^j)
&\leq \max_{s \in [k]} \mathbb{E}(X_i^j \mid X_{i-1}^j, \dots, X_{1}^j, v_i
\in V_s)
= \max_{s \in [k]} \mathbb{E}(X_i^j \mid v_i \in V_s)
\\&= \max_{s \in [k]} \frac{\sum_{x \in F_i}\mathbb{P}(x \in V_j \mid v_i \in
V_s)}{n^{2/3}}
\leq \frac{(1 + \delta/4)|F_i|}{kn^{2/3}}.
\end{align*}
using (b). So, by Lemma~\ref{azuma}, with
probability $1-o(1)$\COMMENT{
So for any fixed $j$, $\mathbb{E}(X_i^j | X_{i-1}^j, \dots, X_{1}^j) \leq a_i =
\frac{(1 + \delta/4)|F_i|}{kn^{2/3}}$, and so $\sum a_i \geq \frac{(1 +
\delta/4)(1-\delta/2k) n}{kn^{2/3}} = \frac{(1 + \delta/4)(1-\delta/2k) m}{n^{2/3}}$. So by
Lemma~\ref{azuma}, $\mathbb{P}(\sum X_i > \frac{(1 + \delta/2) m}{n^{2/3}}) <
\mathbb{P}(\sum X_i > (1+\frac{\delta}{5}) \frac{(1 + \delta/4)(1-\delta/2k) m}{n^{2/3}}) <
e^{-\frac{\delta^2 (1 + \delta/4)(1-\delta/2k) m}{75n^{2/3}}}$ which declines exponentially
with $n$.
} we have that for each~$j$,
$$\sum_{i \in [r]} X_i^j \leq \frac{(1+\delta/2)m}{n^{2/3}}$$
and so for each $j$, the total number of vertices from any of the sets $F_i$
allocated to cluster $V_j$ is at most $(1+\delta/2)m$, proving $(\dagger)$.
\endproof

\subsection{Embedding the vertices of $T$}

Suppose that we have applied the Vertex Allocation Algorithm to find an
approximately uniform allocation of the vertices
of $T$ to the clusters of~$G$. We now wish to embed $T$ in $G$ so that each
vertex is embedded in the cluster to which it is allocated. In principle we
could use the blow--up lemma for this. However numerous complications arise, for
instance because we embed some edges within clusters and because we allow
$\Delta$ to be comparatively large in Section~\ref{sec:unbounded_robust_exp}.
Instead, we embed the vertices of $T$ as follows.

Firstly, to deal with the problem of edges which are allocated within a cluster, 
we shall embed components of $T$ formed by such edges at the same time, using
Theorem~\ref{best_so_far} (it would also be easy to do this directly). To do
this we make the following definition. Let $T$
be a tree on $n$ vertices with root $t_1$, and let the vertices of $T$ be
allocated to clusters $V_1, \dots, V_k$ by a semi-canonical allocation. Then the
\emph{canonical tree} $T_{\textrm{canon}}$ of $T$ is formed by contracting to a 
single vertex each component of the subgraph of $T$ formed of edges which are 
allocated within a cluster. Since the allocation is semi-canonical, each such 
component contains at most $\Delta$ vertices -- we say that these vertices
\emph{correspond} to that contracted vertex in $T_{\textrm{canon}}$. Note also 
that no edge incident to $t_1$ is contracted; we let the root of 
$T_{\textrm{canon}}$ be the vertex corresponding to $t_1$. We shall proceed 
through all of the vertices of $T_{\textrm{canon}}$ in turn using a tidy ancestral 
order, and at time $\tau$ we will embed all of the vertices of $T$ which 
correspond to the vertex $\tau$ of $T_{\textrm{canon}}$ in one step.

Secondly, we must ensure that at each time $\tau$ it is possible to carry out this 
embedding. To do this, each time we embed a vertex $t \in T$ to a vertex $v \in 
G$, we will use Lemma~\ref{sqrt_sets} to select sets $A^+_t$ and $A^-_t$ of 
outneighbours and inneighbours of $v$ in the clusters succeeding and preceding 
that of $v$, each of size at most $2\sqrt{m}$, which are reserved until all of the 
children of $t$ have been embedded. Indeed, while these sets are reserved, no 
vertices may be embedded in them other than 
\begin{itemize}
\item[(i)] the children of $t$, and
\item[(ii)] those vertices of $T$ which correspond to the same vertex of 
$T_{\textrm{canon}}$ as a child of $t$. 
\end{itemize}
We shall refer to these vertices as the \emph{canonical children} of $t$; observe that there are at
most~$\Delta^2$ such vertices. Since we are proceeding
through the vertices of $T_{\textrm{canon}}$ in a tidy ancestral order, this means that at any time 
$\tau$ not too many such sets will be reserved, and so only a small 
proportion of the vertices of any cluster will be reserved. When we later come to 
embed a child~$t'$ of~$t$ for which the edge $tt'$ was allocated canonically, we 
embed $t'$ in $A^+_t$ (if $t \rightarrow t'$) or $A^-_t$ (if $t \leftarrow t'$) in
such a way that we can choose~$A^+_{t'}$ and~$A^-_{t'}$ as desired.

When reading the next algorithm, one should bear in mind that often it is not
apparent that a choice can be made as required by the algorithm. Indeed, if such
a choice is not possible then the algorithm terminates with failure.
Lemma~\ref{embedding_props} will show that under certain conditions on $G$, it will
always be possible to make such choices, and so we can be sure that the
algorithm will succeed.

\medskip
\noindent{ \bf Vertex Embedding Algorithm}

\emph{Input:}
\begin{itemize}
\item A tree $T$ rooted at $t_1$.
\item A constant $\alpha$ and a positive integer $m$.
\item A digraph $G$ on vertex set $V = V_1 \cup \dots \cup V_k$, where each
$V_i$ has size $(1+\alpha)m$, and a
semi-canonical allocation of the vertices of $T$ to the clusters $V_i$, with
$t_1$ allocated to~$V_1$.
\item Finally, a vertex $v^* \in V_1$ to which $t_1$ should be embedded, and
constants $c$ and $\gamma$.
\end{itemize}

\emph{Initialisation:} Form the canonical tree $T_{\textrm{canon}}$ of $T$ as explained
above, and choose a tidy ancestral ordering $1, 2, \dots, n'$ of the vertices of
$T_{\textrm{canon}}$. Let $t_1, \dots, t_n$ be a corresponding order of the vertices of
$T$ (so if $t_i \in T$ corresponds to $i \in T_{\textrm{canon}}$ and $t_j \in T$
corresponds to $j \in T_{\textrm{canon}}$ then $t_i$ appears before $t_j$ if and only if
$i<j$.)

\emph{Procedure:} At time $\tau$ we shall embed the vertices $t_r, \dots,
t_{r+s-1}$ of $T$ corresponding to vertex~$\tau$ of $T_{\textrm{canon}}$. Each vertex $t_i$
will
be embedded to a vertex $v_i$ of $G$, where $v_1 = v^*$. Then, for each~$t_i$ we
will reserve sets
$A_{t_i}^+$ and $A_{t_i}^-$ of vertices of $G$ for the canonical children of $t_i$. To do
this, at each time $\tau$ with $ 1 \leq \tau \leq n'$, take the following steps.
\begin{itemize}
\item[(1)] We say that a vertex $t_i$ of $T$ is \emph{open} at time $\tau$ if $t_i$ has
been embedded but some child of $t_i$ has not yet been embedded. Define the set
$B^\tau$ of vertices of $G$ unavailable for use at time $\tau$ to consist of the
vertices already occupied and the sets reserved for the canonical children of open
vertices, so
$$ B^\tau = \{v_1, \dots, v_{r-1}\} \cup \bigcup_{t_i \colon t_i \textrm{ is
open}} (A^+_{t_i} \cup A^-_{t_i}).$$
For each cluster $V_j$, let  $V_j^\tau =V_j \sm B^\tau$, so $V_j^\tau$ is the
set of available vertices of $V_j$.
\item[(2)] If $\tau=1$ embed $t_1$ to $v_1$. Alternatively, if $\tau > 1$:
\begin{itemize}
\item[(2.1)] Precisely one of the vertices $t_{r}, \dots, t_{r+s-1}$ of $T$
corresponding to vertex $\tau$ of $T_{\textrm{canon}}$ has a parent already embedded; we
may assume this vertex is $t_r$. Let $t_p$ be the already-embedded parent
(so $p < r$, and when $t_p$ was embedded sets $A_{t_p}^+$ and $A_{t_p}^-$ were
chosen). Let $V_j$ be the cluster to which $t_p$ is embedded.
\item[(2.2)] If $t_p \rightarrow t_{r}$, choose a set $S$ of $3s$ vertices of
$A^+_{t_p}\subseteq V_{j+1}$ such that for each $v \in S$
$$|N^+(v) \cap V_{j+2}^\tau| \geq \gamma m \textrm{ and }|N^-(v) \cap
V_{j}^\tau| \geq \gamma m.$$
If $t_p \leftarrow t_{r}$, choose a set $S$ of $3s$ vertices of
$A^-_{t_p}\subseteq V_{j-1}$ so for each $v \in S$
$$|N^+(v) \cap V_{j}^\tau| \geq \gamma m \textrm{ and }|N^-(v) \cap
V_{j-2}^\tau| \geq \gamma m.$$
\item[(2.3)] Then choose a copy of $T[t_r, \dots, t_{r+s-1}]$ in $G[S]$, and embed
each vertex $t_{i}$ to the corresponding vertex $v_i$ in this copy.
\end{itemize}
\item[(3)] In step (2), we embedded each of $t_r, \dots, t_{r+s-1}$ in the same
cluster; let $V_q$ be this cluster. For each $r \leq i \leq r+s-1$, choose
sets
$$A^+_{t_i} \subseteq N^+(v_i) \cap V_{q+1}^\tau \text{ and } A^-_{t_i}
\subseteq N^-(v_i) \cap V_{q-1}^\tau$$ such that the sets $A^+_{t_i}$ and
$A^-_{t_i}$ are all pairwise disjoint, each $A^+_{t_i}$ and each $A^-_{t_i}$ is
$(c, \gamma)$-good, and $|A^+_{t_i}|, |A^-_{t_i}| \leq
2\sqrt{m}$ for each $i$.
\end{itemize}
Whenever there are several choices (for example if there are several
possibilities for $S$ in (2.2)), take the lexicographically first of these. This
ensures that for each input, the output is uniquely defined (i.e. we can view
the algorithm as being deterministic).

\emph{Termination:} If at any point it is not possible to make the choice
required, terminate with failure. Otherwise, terminate after every vertex of
$T_{\textrm{canon}}$ has been processed, at which point $\psi(t_i) = v_i$ for each
$t_i \in
T$ is an embedding $\psi$ of $T$ into $G$, by construction.

\begin{lemma} \label{embedding_props}
Suppose that $1/n \ll 1/\Delta, 1/k \ll \eps \ll \gamma \ll c \ll d \ll \alpha \leq 2$, and
let $m = n/k$.\begin{itemize}
\item[(1)] Let $T$ be a directed tree on at most $n$ vertices with root $t_1$
and $\Delta(T)
\leq \Delta$.
\item[(2)] Let $G$ be an $\eps$-regular $d$-dense cycle of cluster tournaments
on clusters $V_1, \dots, V_k$, each of size $(1+\alpha)m$, and let $v^* \in V_1$
have at least $\gamma m$ inneighbours in $V_k$ and at least $\gamma m$
outneighbours in $V_2$.
\item[(3)] Let the vertices of $T$ be allocated to the clusters $V_1, \dots,
V_k$ so that at most $(1+\alpha/2)m$ vertices are allocated to any one cluster
$V_i$, and so that the allocation is semi-canonical.
\end{itemize}
Then if we apply the Vertex Embedding Algorithm to $T$ and $G$ with this
allocation and constants $c$ and $\gamma$, then it will successfully embed $T$
into $G$ with $t_1$ embedded to $v^*$.
\end{lemma}

\proof The Vertex Embedding Algorithm will
only fail if at some point it is not possible to make
the required choice. So to demonstrate that the algorithm will succeed, it is
enough to show that it is always possible to make the required choices.

In the initialisation we are required to choose a tidy ancestral ordering of the vertices
of the rooted tree $T_{\textrm{canon}}$; the existence of such a choice is guaranteed by
Lemma~\ref{open_vertices}. Now, consider the set of unavailable vertices
$B^\tau$ at some time $\tau$. Since the Vertex Embedding Algorithm embeds each
vertex in the cluster to which it was allocated, we know that at most
$(1+\alpha/2)m$ vertices of each $V_j$ are already occupied. Furthermore,
suppose that vertex $t_i$ of~$T$ is open at time $\tau$. Then $t_i$ must
correspond to a vertex $\tau' < \tau$ of $T_{\textrm{canon}}$, such that $\tau'$ has a child
$\tau'' \geq \tau$. Since we are processing the vertices of $T_{\textrm{canon}}$ in a tidy
order, there can be at most $\log_2 n' \leq \log_2 n$ such vertices of $T_{\textrm{canon}}$. As
each vertex of $T_{\textrm{canon}}$ corresponds to at most $\Delta$ vertices of $T$, at most
$\Delta \log_2 n$ vertices of $T$ are open at time $\tau$. Therefore, at any
time $\tau$, the total number of vertices in reserved sets $A_{t_i}^+$ and
$A_{t_i}^-$ is at most $4 \Delta \sqrt{m} \log_2 n  \leq \alpha m/4$. So for any
cluster $V_j$, at any time $\tau$ at most $(1+\alpha/2)m + \alpha m/4$ vertices
of $V_j$ are unavailable, and so $|V_j^\tau| \geq \alpha m/4$.

We can now demonstrate that it is possible to make the other choices that the
algorithm asks for. Indeed, in step (2.2), if $t_p \rightarrow t_{r}$ with $t_p$
embedded into $V_j$, then the
algorithm has to choose a set $S$ of $3s \leq 3\Delta$ vertices of
$A^+_{t_p}$ such that each $v \in S$ has $|N^+(v) \cap V_{j+2}^\tau| \geq \gamma
m$ and $|N^-(v) \cap V_{j}^\tau| \geq \gamma m$. But when $A^+_{t_p}$ was chosen
at an earlier time $\tau'$, it was chosen to be $(c, \gamma)$-good. Since the
vertex $v_p$ to which $t_p$ was embedded is in cluster $V_j$, $A^+_{t_p}
\subseteq V_{j+1}$. Moreover, since $|V_j^\tau| \geq \alpha m/4 \geq (1+\alpha)cm$ and
$|V_{j+2}^\tau| \geq \alpha m/4 \geq (1+\alpha)cm$, $A^+_{t_p}$ must contain at least
$\gamma \sqrt{m}$ vertices $v$ such that $|N^+(v) \cap V_{j+2}^\tau| \geq \gamma
m$ and $|N^-(v) \cap V_{j}^\tau| \geq \gamma m$. Furthermore, since $t_r$ is a
child of $t_p$, $t_p$ has been open since its embedding, and so only canonical children of
$t_p$ (of which there are at most $\Delta^2$) can have been embedded in
$A^+_{t_p}$. So it is indeed possible to select such a set $S$ of $3s$ vertices
as required. The argument for the case when $t_p \leftarrow t_{r}$ is similar.
As for (2.3), observe that $G[S]$
is a tournament on $3s$ vertices, and that $T[t_r, \dots, t_{r+s-1}]$
is a directed tree on $s$ vertices. So by Theorem~\ref{best_so_far}, $G[S]$
contains a
copy of $T[t_r, \dots, t_{r+s-1}]$, so we may choose such a copy.

Finally we come to step (3). In this step we have just embedded at most $\Delta$
vertices $t_{r}, \dots, t_{r+s-1}$ in some cluster $V_q$, and we wish to choose
sets $A_{t_i}^+$ and $A_{t_i}^-$ for each such vertex~$t_i$. When embedding
these vertices we ensured that for each $i$ the vertex $v_{i}$ to which $t_i$
was embedded satisfied $|N^+(v_i) \cap V_{q+1}^\tau| \geq \gamma m$ (for
$\tau = 1$ this holds instead by the condition on the outneighbours of $v^* =
v_1$). So suppose we have chosen $A_{t_r}^+, A_{t_{r+1}}^+, \dots,
A_{t_{r+\ell-1}}^+$ and we now wish to choose $A_{t_\ell}^+$. Then the
previously chosen $A_{t_i}^+$ contain at most $2\Delta\sqrt{m}$ vertices between
them, and so at least $3\gamma m/4 \geq (1+\alpha)\gamma m/2$ vertices of
$N^+(v_\ell) \cap V_{q+1}^\tau$ have not been used in these previous sets. So by
Lemma~\ref{sqrt_sets}, we may choose a $(c, \gamma)$-good set $A_{t_\ell}^+
\subseteq N^+(v_\ell) \cap V_{q+1}^\tau$ of size at most $2\sqrt{m}$ which is
disjoint from all of the previously chosen $A_{t_i}^+$. Do this for each vertex
$t_i$ in turn; the choice of the sets $A_{t_i}^-$ is similar.
\endproof

We can now give the proof of the main lemmas of this section, beginning with the
proof of Lemma~\ref{cluster_cycle_bounded_deg}.

\medskip\noindent {\bf Proof of Lemma~\ref{cluster_cycle_bounded_deg}.}
Apply the Vertex Allocation Algorithm to allocate the vertices of~$T$ to the
clusters $V_1, \dots, V_k$. Then by Lemma~\ref{allocation_props}(a) this
allocation is semi-canonical, and by Lemma~\ref{allocation_props}(c) at most
$(1 + \alpha/2)m$ vertices are allocated to each of the
$k$ clusters $V_i$. Next, apply the Vertex Embedding Algorithm to $T$ and $G$,
giving this allocation as input. By Lemma~\ref{embedding_props}, this will
successfully embed $T$ in $G$ with $t_1$ embedded to $v^*$.
\endproof

\medskip\noindent {\bf Proof of Lemma~\ref{robust_expander_case_bounded_deg}.}
If $\alpha >2$ then $G$ contains a copy of $T$ by Theorem~\ref{best_so_far}. So
we may assume that $\alpha \leq 2$. We begin by introducing new constants $1/n
\ll 1/M \ll 1/M' \ll \eps \ll \eps'\ll d \ll
\mu$. Then by Lemma~\ref{KOT_combined}, $G$ contains an $\eps$-regular $d$-dense
cycle of cluster tournaments~$G'$ on clusters $V_1, \dots, V_k$, where $M' \leq
k \leq M$, and $|V_1| = \dots = |V_k| \geq (1-\eps)(1+\alpha)n/k \geq
(1+\alpha/2)n/k$. For each $i$ choose $V_i' \subseteq V_i$ of size $|V_i'| =
(1+\alpha/2)n/k$ uniformly at random. By Lemma~\ref{regular_restriction} we may
fix an outcome of these choices so that $G'' = G'[V_1' \cup \dots \cup V_k']$ is
a $\eps'$-regular $d/2$-dense cycle of cluster tournaments. So by
Lemma~\ref{cluster_cycle_bounded_deg} $G''$ contains a copy
of~$T$, so $G$ contains $T$ also.
\endproof

We finish this section with an analogous result to
Lemma~\ref{cluster_cycle_bounded_deg} for small trees (i.e. the result does not
demand that $|T|$ is large compared to $|G|$).

\begin{lemma} \label{embed_small_trees}
Suppose that $1/m \ll 1/k, 1/\Delta \ll \eps \ll d \ll \alpha \leq 2$, and that $1/k
\ll \delta$. Let $G$ be an $\eps$-regular $d$-dense cycle of cluster tournaments
on clusters $V_1, \dots, V_k$, each of size $(1+\alpha) m$, and let $v^* \in
V_1$ have at least $d^2m$ inneighbours in $V_k$ and at least $d^2m$
outneighbours in $V_2$. Let $T$ be a directed tree on at most $m$ vertices,
rooted at $t_1$ and with $\Delta(T) \leq \Delta$, and let $T^{\textrm{far}}$ be the
subgraph of $T$ induced by the vertices $x \in T$ with $d(t_1, x) \geq k^3$. Let
$\G_T$ denote the set of copies of $T$ in $G$ for which the vertex $t_1$ of $T$
corresponds to vertex $v^*$ of $G$. Then $\G_T$ is non-empty. Furthermore, there
exists a probability distribution on $\G_T$ such that if a member of $\G_T$ is
selected at random according to this distribution, then for each $i$,
$$\mathbb{E}(\# \textrm{ vertices of }T^{\textrm{far}} \textrm{ embedded in }V_i) \leq
\frac{(1 + \delta)|T^{\textrm{far}}|}{k}.$$
\end{lemma}

The probability distribution will actually be, for each member of $\G_T$, the
probability that applying first the Vertex Allocation Algorithm and then the
Vertex Embedding Algorithm gives this copy of $T$ in $G$ (recall that actually the
Vertex Embedding Algorithm is purely deterministic).

\proof
Apply the Vertex Allocation Algorithm to allocate the vertices of $T$ to
the clusters $V_i$. Since $|T| \leq m$, at most $m$ vertices can be allocated to any
cluster, and the allocation obtained is semi-canonical by
Lemma~\ref{allocation_props}(a). Next, introduce constants $\eps \ll \gamma \ll
c \ll d$, and apply the Vertex Embedding Algorithm to embed $T$ in $G$. By
Lemma~\ref{embedding_props}, this will successfully embed $T$ in $G$, with $t_1$
embedded to $v^*$, and every vertex of $T$ embedded in the cluster to which it
was allocated. So it remains only to show that for each $i$, the expected number
of vertices of $T^{\textrm{far}}$ allocated to $V_i$ is at most $(1+\delta)|T^{\textrm{far}}|/k$.
But since for any $x \in T^{\textrm{far}}$ we have $d(x, t_1) \geq k^3$, by
Lemma~\ref{allocation_props}(b) applied with $u = t_1$,
$$ \mathbb{P}(x\textrm{ is allocated to } V_i) = \frac{(1\pm\delta)}{k}$$
for each $i$, and the result follows.
\endproof

\section{Embedding trees of unbounded maximum degree in a robust
outexpander.}\label{sec:unbounded_robust_exp}

\subsection{Section outline}

Having proved the desired result for trees of bounded maximum degree, we now
move onto proving a similar result for trees with no such bound, with a constant
of $2$ rather than $1$ in the condition on the order of $G$. This is the
following lemma.

\begin{lemma} \label{robust_expander_case}
Suppose that $1/n \ll \mu \ll \nu \ll \eta \ll \alpha$, that $G$ is a
tournament on $2(1+\alpha)n$ vertices which is a robust $(\mu, \nu)$-outexpander
with $\delta^0(G) \geq \eta n$ and that $T$ is a directed tree on~$n$ vertices.
Then $G$ contains a copy of $T$.
\end{lemma}

To prove this, we shall begin with a definition. In
Section~\ref{subsec:core_tree} we shall define the core tree~$T_c$ of a tree
$T$. This is a subtree of $T$ which has bounded maximum degree, and the property
that all components of $T - T_c$ are small. Then in 
Section~\ref{subsec:extendedtree} we will show that $T_c$ can be extended to an 
`extended tree' $T_{\textrm{ext}}$ which also has bounded maximum degree, and also has the 
property that few vertices of $T_{\textrm{ext}}$ have neighbours outside $T_{\textrm{ext}}$. We will 
embed the extended tree $T_{\textrm{ext}}$ by a similar method to that of the previous 
section. We will need to do this so that the small number of vertices of 
$T_{\textrm{ext}}$ with neighbours outside $T_{\textrm{ext}}$ are embedded to vertices of $G$ 
with large in- and outdegree in $G$. 
In Section~\ref{subsec:rob_exp_proof1} we will use
our results from Section~\ref{sec:bounded_robust_exp} to prove
Lemma~\ref{embed_with_restrictions} on embedding trees of bounded maximum
degree. This is similar to Lemma~\ref{cluster_cycle_bounded_deg}, but allows us
also to demand that a small subset $H \subseteq V(T)$ of the vertices of $T$,
satisfying certain conditions, should be embedded in a small subset $U$ of the
vertices of $G$. This will allow us to embed $T_{\textrm{ext}}$ in $G$ in the desired manner. Finally, in Section~\ref{subsec:rob_exp_proof2} we will
complete the proof of Lemma~\ref{robust_expander_case} by first using
Lemma~\ref{embed_with_restrictions} to embed $T_{\textrm{ext}}$ into~$G$ as described and then embedding each component of $T-T_{\textrm{ext}}$ in 
the
unoccupied vertices of $G$.

\subsection{The core tree} \label{subsec:core_tree}

Let $T$ be a tree on $n$ vertices, and let $\Delta \geq 2$ be fixed. Then we say
that a vertex $v$ of $T$ is \emph{$\Delta$-core} if every edge incident to $v$
has $w_e(v) \leq (1-1/\Delta)n$. We call the subgraph of $T$ induced by
$\Delta$-core vertices of $T$ the \emph{core tree of $T$ with parameter
$\Delta$}, and denote it by $T_c$ (the value of $\Delta$ will always be clear
from the context). With this definition, for any tree $T$, the core tree of $T$
with parameter $\Delta$ is the same as the $\Delta$-heart of $T$ considered by
Thomason and H\"aggkvist in~\cite{HT}. The statements of the next proposition are also noted in Section~3 of~\cite{HT}, but we include the proof for completeness. 

\begin{prop} \label{core_tree_props}
Let $T$ be a tree on $n$ vertices, let $\Delta \geq 2$ and let $T_c$ be the core
tree of $T$ with parameter $\Delta$. Then:
\begin{itemize}
\item[(i)] $T_c$ is a tree containing at least one vertex.
\item[(ii)] $w_e(x) \geq n/\Delta$ if $e = xy$ is an edge of $T_c$.
\item[(iii)] $\Delta(T_c) \leq \Delta$.
\item[(iv)] Every component subtree $T'$ of $T-T_c$ has $|T'| \leq n/\Delta$.
\end{itemize}
\end{prop}

\proof For (i), note that since $\Delta \geq 2$, for any edge $e=uv$ of $T$ at
most one of $w_e(u) > (1-1/\Delta)n$ and $w_e(v) > (1-1/\Delta)n$ holds. Since
$T$ has more vertices than edges, there must therefore be some vertex $v \in T$
such that $w_e(v) \leq (1-1/\Delta)n$ for every edge $e$ incident to $v$, and so
$v \in T_c$. It remains to show that $T_c$ is connected. Observe that if $u, v,
w$ are distinct vertices of $T$ such that there is an edge between $u$
and $v$ and an edge between $v$ and $w$, then $w_{uv}(u) > w_{vw}(v)$. Now,
suppose $x, y \in T_c$, and let $x=v_1, v_2, \dots, v_r=y$ be the vertices of
the path from $x$ to $y$ in $T$ (in order). Suppose for a contradiction that
some $v_i$ is
not in $T_c$. Then for some neighbour $z$ of $v_i$, $w_{v_iz}(v_i) >
(1-1/\Delta)n$. If $z \neq v_{i+1}$, then for each $i \leq j \leq r-1$ we have
$w_{v_jv_{j+1}}(v_{j+1}) > (1-1/\Delta)n$, and so $y \notin T_c$, giving a
contradiction. On the
other hand, if $z = v_{i+1}$, then for each $2 \leq j \leq i$,
$w_{v_{j-1}v_{j}}(v_{j-1}) > (1-1/\Delta)n$, and so $x \notin T_c$, again giving
a contradiction.

Now, (ii) is immediate from the fact that if $e=xy$ is an edge of $T$ then
$w_e(x) + w_e(y) = n$. Then (iii) follows directly from (ii), as the sum of
$w_e(v)$ over all edges incident to $v$ is $n-1$.

Finally, for (iv), observe that for any such $T'$ there is $u \in T'$, $v \in
T_c$ with $e = uv$ an edge of $T$. Suppose that $|T'| > n/\Delta$. Then $w_e(v)
\geq |T'| > n/\Delta$, and so $w_e(u) \leq (1-1/\Delta)n$. But since $w_{e'}(u)
< w_e(v) \leq (1-1/\Delta)n$ for every other edge $e'$ incident to $u$, this
means that $u \in T_c$, giving a contradiction.
\endproof

Note that $T_c$ is an undirected tree obtained from an undirected tree $T$.
However we will often refer to the core tree of a directed tree $T$; this means
the directed tree formed by taking the core tree $T_c$ of the underlying graph
$T_{\textrm{under}}$ (an undirected tree) and directing each edge of~$T_c$ as it is directed in $T$.

The idea behind this definition is that the core tree is a bounded degree tree. The
general technique we shall use to work with a tree $T$ of unbounded maximum
degree (in both this and later sections) is to first consider the core tree
$T_c$, and then consider separately each component of $T - T_c$, making use of
the fact that each such component is small.

\subsection{Leading paths} \label{subsec:leading_paths}

Let $T$ be a tree on $n$ vertices, rooted at $t_1$, let $H \subseteq V(T)$, and
let $k$ be a positive integer. For any vertex $x \in T$, there is a unique path
in $T$ from $x$ to $t_1$; let~$P_x$ denote the set of the first $k$ vertices of
this path, starting from $x$. Let $H^1 = \bigcup_{x \in H} P_x$, and then for
each $i \geq 1$ let $H^{i+1}$ be formed from $H^{i}$ by adding the vertices of
$P_x$ for any $x \in H^{i}$ with at least two children in $H^{i}$. After at most
$n$ steps we must have $H^{i} = H^{i+1}$, when we terminate the process. We
refer to this final $H^{i}$ as \emph{$H$ with leading paths included},
denoted~$\mathcal{P}_k(H)$. So $H \subseteq \mathcal{P}_k(H) \subseteq V(T)$.
Note that $\mathcal{P}_k(H)$ depends on both the value of $k$ and the root $t_1$
of $T$.

Next we shall prove two results which will enable us to make use of this
definition. The first shows that if $H$ is small then $\mathcal{P}_k(H)$ is
small, and the second shows that if $H$ is small then it is possible to embed
any component $T'$ of $T[\mathcal{P}_k(H)]$ in a regular and dense cycle of
cluster tournaments such that the vertices of $V(T') \cap H$ are embedded in the
first cluster and the `root' of $H$ is embedded in a given cluster.

\begin{prop} \label{size_of_p(h)}
Let $k$ be any positive integer, let $T$ be a tree on $n$ vertices, rooted at
some $t_1 \in T$, and let $H \subseteq V(T)$. Then $|\mathcal{P}_k(H)| \leq
3k|H|$.
\end{prop}

\proof Consider any component $T'$ of $T[\mathcal{P}_k(H)]$, and let $t_1'$ be
the unique vertex of $T'$ with minimal $d(t_1, t_1')$. Then every vertex of $T'$
lies on the path from some vertex of $H$ to $t_1$, and so~$T'$ is precisely the
set of vertices in paths between $t_1'$ and vertices of $H \cap V(T')$. Thus
only~$t_1'$ and vertices of $H$ can be leaves of $T'$. It follows that
$T[\mathcal{P}_k(H)]$ has at most $2|H|$ leaves. Since $T[\mathcal{P}_k(H)]$ is
a forest, it follows that the number of vertices of $T[\mathcal{P}_k(H)]$ with at least two
children in $T[\mathcal{P}_k(H)]$ is also at most $2|H|$. Furthermore, any
vertex $x \in T$ for which the vertices of $P_x$ were added to
$\mathcal{P}_k(H)$ at any stage is either a member of $H$ or has at least two
children in $\mathcal{P}_k(H)$. This is true for at most $3|H|$ vertices $x$,
and for each such vertex at most $k$ vertices were added.
\endproof

\begin{lemma} \label{embed_p(h)}
Suppose that $1/m \ll 1/k \ll \eps \ll d$. Let $T$ be a directed tree rooted at
some $t_1 \in T$. Let $H \subseteq V(T)$ be of size $|H| \leq m/10k$, let $T'$
be a component of $T[\mathcal{P}_k(H)]$ which does not contain~$t_1$, and let
$t'_1$ be the unique vertex of $T'$ with minimal $d(t'_1, t_1)$. Let $G$ be an
$\eps$-regular $d$-dense cycle of cluster tournaments on clusters $V_1, \dots,
V_k$, each of size $m$. Then for any $j \in [k]$, $G$ contains a copy of $T'$
with the vertex $t'_1$ corresponding to some vertex of~$V_j$, and every vertex
in $V(T') \cap H$ corresponding to some vertex of $V_1$.
\end{lemma}

\proof
Informally, from the perspective of $t'_1$, $T'$ begins with a path of length
$k-1$ (from $t'_1$ to $t$, say) before possibly branching out. So we shall find
a copy of $T'$ in $G$ by first embedding the vertices of this path so that
$t'_1$ is embedded in $V_j$ and $t$ is embedded in $V_1$. We then embed all of
the remaining vertices of $T'$ in $V_1$.

More formally, note that for each $0 \leq s \leq k-1$ there is precisely one
vertex $x_s$ of $T'$ with $d(t'_1, x_s) = s$ (so $x_0 = t'_1$, and $x_i \notin H$
for any $i < k-1$). Let $F \subseteq [k-1]$ be the set of those $s$ such that
$x_{s-1} \rightarrow x_s$, and let $B \subseteq [k-1]$ be the set of $s$ such
that $x_{s-1} \leftarrow x_s$. Then $|F|+|B| = k-1$, so either $|F| > k-j$ or
$|B| \geq j-1$. Suppose first that $|B| \geq j-1$. Then choose $B' \subseteq B$
of size $j-1$. We shall allocate the vertices of $T'$ to the clusters $V_1,
\dots, V_j$. Begin by allocating $x_0$ to $V_j$. Then for each $s \in [k-1]$ in
turn, let $V_i$ be the cluster to which $x_{s-1}$ was allocated, and allocate
$x_s$ to $V_i$ if $s \notin B'$, or to $V_{i-1}$ if $s \in B'$. Then since $|B'|
= j-1$, $x_{k-1}$ will be assigned to $V_1$. Finally, allocate all other
vertices of $T'$ to $V_1$. Then every edge of $T'$ is allocated either
canonically or within a cluster.

Next we shall embed $T'$ in $G$ so that every vertex is embedded within the
cluster to which it is allocated. To begin, by a standard regularity argument we
may choose for each $i$ a set $V_i' \subseteq V_i$ so that $|V_i'| \geq 9m/10$
and every vertex $v \in V_i'$ has at least $dm/2$ outneighbours in $V_{i+1}'$.
Let $G' = G[V_1' \cup \dots \cup V_k']$. Now, for each $i$, let $S_i$ be the set
of vertices of $T'$ allocated to $V_i$. So $|S_2|, \dots, |S_k| \leq k-1$ and $|S_1| \leq |T'|$. Then by Proposition~\ref{size_of_p(h)},
$3|S_1| \leq 3|T'| \leq 9k|H| \leq |V_1'|$. So by Theorem~\ref{best_so_far} we
may embed $T'[S_1]$ in $G[V_1']$. Now suppose that we have successfully embedded
$T'[S_1 \cup \dots \cup S_{i-1}]$ in $G[V_1' \cup \dots \cup V_{i-1}']$ for some
$i \leq j$. Then precisely one vertex $t \in S_i$ has a neighbour $t' \in
S_{i-1}$, and $t'$ has already been embedded to some $v' \in V_{i-1}'$. Now $v'$
has at least $dm/2 \geq 3|S_i|$ outneighbours in $V_i'$, and so by
Theorem~\ref{best_so_far} we may embed $T'[S_i]$ among these outneighbours. Let
$v$ be the vertex to which $t$ is embedded; then since $v$ is an outneighbour of
$v'$, we have extended our embedding to an embedding of $T'[S_1 \cup \dots \cup
S_i]$ in $G[V_1' \cup \dots \cup V_i']$. Continuing in this manner we obtain an
embedding of $T'$ in $G$, with $t'_1$ embedded in $V_j$ and $V(T') \sm \{x_0,
\dots, x_{k-2}\} \supseteq V(T') \cap H$ embedded into $V_1$, as desired. A
similar argument will achieve this if $|F| > k-j$.
\endproof

\subsection{The extended tree.} \label{subsec:extendedtree}
The next lemma combines the ideas of the core tree and leading paths to give the 
structure within a tree $T$ which we shall use to prove 
Lemma~\ref{robust_expander_case}. It shows that given a tree $T$ we may extend the 
core tree $T_c$ of $T$ with parameter $\Delta$ to an `extended tree' $T_{\textrm{ext}}$ 
which, like $T_c$, has bounded maximum degree (although this bound is now much 
larger than $\Delta$). $T_{\textrm{ext}}$ will also have the property that only a small 
subset $H$ of the vertices of $T_{\textrm{ext}}$ have neighbours outside $T_{\textrm{ext}}$, and 
that few vertices of $T_{\textrm{ext}}$ are close to a vertex of $\mathcal{P}_k(H)$.

\begin{lemma} \label{extendedtree}
Suppose that $1/n, 1/\Delta^* \ll 1/\Delta, 1/k, \omega \ll 1$. Let $T$ be a tree 
on $n$ vertices, and let $T_c$ be the core tree of $T$ with parameter $\Delta$. 
Choose any vertex $t_1 \in T_c$ as the root of~$T$. Then there exists a subtree 
$T_{\textrm{ext}}$ of $T$ and a subset $H \subseteq V(T_{\textrm{ext}})$ which satisfy the following 
properties.
\begin{itemize}
\item[(i)] $T_c \subseteq T_{\textrm{ext}}$.  
\item[(ii)]  $\Delta(T_{\textrm{ext}}) \leq \Delta^*$.
\item[(iii)]  For any edge $e$ between $V(T - T_{\textrm{ext}})$ and $V(T_{\textrm{ext}})$, the endvertex of $e$ in $V(T_{\textrm{ext}})$ lies in~$H$.
\item[(iv)]  The number of
vertices $v \in T_{\textrm{ext}}$ which satisfy $1 \leq d(v, \mathcal{P}_k(H)) \leq
k^3$ is at most $\omega n$.
\item[(v)]  $|H| \leq n/\Delta^{k^{1/\omega}}$.
\end{itemize}
\end{lemma}

\proof We consider the subgraph $T-E(T_c)$ of $T$ obtained by deleting the edges 
(but not the vertices) of $T_c$ from~$T$. Each vertex $v \in T_c$ lies in a 
separate component of $T-E(T_c)$; we denote the component containing $v$ by $T_v$. 
Then the trees $T_c$ and $\{T_v: v \in T_c\}$ partition the edges of $T$, and the 
trees $\{T_v: v \in T_c\}$ partition the vertices of $T$.

We say that a vertex $v \in T_c$ is \emph{$i$-heavy} if $|T_v| \geq \Delta_i := 
\Delta^{k^i}$. For any integer~$i$, let $H_i$ denote the set of $i$-heavy vertices 
in $T_c$. So\COMMENT{
Since the trees $T_v: v \in H_i$ are pairwise vertex-disjoint and each such tree 
contains at least
$\Delta_i$ vertices.
} $|H_i| \leq n/\Delta_i$, and so by Proposition~\ref{size_of_p(h)} we have $|\mathcal{P}_k(H_i)| \leq 3kn/\Delta_i$ for each $i$.
We wish to choose a large integer $t$ so that few vertices of $T$ lie in trees $T_v$ for which $v$ is not in $H_t$ but is close to a member of $\mathcal{P}_k(H_t)$. The next claim shows that this is possible.

\medskip
\noindent
{\bf Claim.}
For some natural number $1/\omega \leq t \leq 3/\omega$ we have
\begin{equation} \label{eq:claim}
\big| \bigcup_{\substack{v \in V(T_c) \sm H_t \\ d(v, \mathcal{P}_k(H_t)) \leq k^3}} T_v \big| \leq \omega n.
\end{equation}
\medskip

\noindent
{\bf Proof of Claim.}
Observe that for each integer $i$ with $1/\omega \leq i \leq 3/\omega$, if $v \in V(T_c) \sm
H_{i-1}$ then
$|T_v| < \Delta_{i-1}$, and so\COMMENT{Using $\Delta^{k^2} \geq 450\times 32 
k^2/\alpha \lambda d = 450k/\alpha\delta$. Also that $k \geq 2, i \geq 2$ in the 
previous step.}
$$\big|\bigcup_{\substack{ v \in V(T_c) \sm H_{i-1} \\ d(v, \mathcal{P}_k(H_i)) 
\leq k^3}} T_v\big| < |\mathcal{P}_k(H_i)| \Delta^{k^3+1}\Delta_{i-1} \leq
\frac{3k\Delta^{k^3+1}\Delta^{k^{i-1}}n}{\Delta^{k^{i}}}\leq
\frac{3kn}{\Delta^{k^i/2}} \leq \omega n/3.$$
Now let
$$B_i := \bigcup_{\substack{v \in H_{i-1} \sm H_i \\ d(v, \mathcal{P}_k(H_i)) \leq 
k^3}} T_v.$$
Then the sets $B_i$ are pairwise disjoint\COMMENT{Since $H_1 \supseteq H_2 
\supseteq \dots$.} subsets of $V(T)$. If the claim is false, then $|B_i| > 
2\omega n/3$ for every integer $i$ with $1/\omega \leq i \leq 3/\omega$, and so
$|\bigcup_{1/\omega \leq i \leq 3/\omega} B_i| > n$, giving a contradiction.
\endproof

Fix such a value of $t$, and let $H = H_t$. We define the extended tree
$T_{\textrm{ext}}$ by $T_{\textrm{ext}} := T_c \cup \bigcup_{v \in V(T_c) \sm H} T_v$. Then 
$T_{\textrm{ext}}$ is a subtree of $T$ with $T_c \subseteq T_{\textrm{ext}}$, so (i) is satisfied. 
Since $H \subseteq V(T_c)$, we have $H \subseteq V(T_{\textrm{ext}})$ as desired. Also (ii) 
holds since any vertex $u \in T_{\textrm{ext}}$ has at most $\Delta$ neighbours in $T_c$ 
and at most $\Delta_t$ neighbours in the single tree $T_v$ with $v \in T_c$ which 
contains $u$. So $\Delta(T_{\textrm{ext}}) \leq \Delta + \Delta_t \leq \Delta + 
\Delta^{k^{3/\omega}} \leq \Delta^*$. For (iii), observe that if $u \notin 
T_{\textrm{ext}}$, then $u$ must lie in some $T_v$ with $v \in H$. But then if $u$ has a 
neighbour in $T_{\textrm{ext}}$ this neighbour must be $v$. For (iv), consider any $u \in 
T_{\textrm{ext}}$ satisfying $1 \leq d(u, \mathcal{P}_k(H)) \leq k^3$. Since $d(u, 
\mathcal{P}_k(H)) \geq 1$ we know that $u \notin H_t$, so if $u \in T_c$ then $u$ 
is counted in (\ref{eq:claim}). If $u \notin T_c$ then there exists $v$ such that 
$u \in T_v$ and $v \in V(T_c) \sm H$. Note that $\mathcal{P}_k(H)\subseteq V(T_c)$ 
(since $t_1\in T_c$). This in turn implies that $d(v, \mathcal{P}_k(H))<d(u, 
\mathcal{P}_k(H)) \leq k^3$. So $u$ is also counted in (\ref{eq:claim}). Finally, 
for (v), recall that $|H| \leq n/\Delta_t \leq n/\Delta^{k^{1/\omega}}$.
\endproof

\subsection{Embedding trees of bounded maximum degree with restrictions}
\label{subsec:rob_exp_proof1}

In this section we shall prove the following lemma, which is similar to
Lemma~\ref{cluster_cycle_bounded_deg}, but which allows us to restrict some
vertices of $T$ to a subset of $V(G)$.

\begin{lemma} \label{embed_with_restrictions}
Suppose that $1/n \ll 1/\Delta, 1/k \ll \eps \ll d \ll \alpha, \lambda
\leq 1/2$, that $m = n/k$, that $\lambda \leq \alpha/4$ and that $\delta =
d\lambda/8k$. Let $T$ be a directed tree on $n$
vertices rooted at $t_1$ and with $\Delta(T) \leq \Delta$. Let $H \subseteq
V(T)$ be such that  $|H| \leq \delta n/7k$ and $|\{x \in T \colon 1 \leq d(x,
\mathcal{P}_k(H)) \leq k^3\}| \leq \delta n$. Let $G$ be an $\eps$-regular
$d$-dense cycle of cluster tournaments on clusters $V_1, \dots, V_k$, each of
size $(1+\alpha)m$, and let $U \subseteq V_1 \cup \dots \cup V_k$ have size $|U|
\geq \lambda n$. Then $T$ can be embedded in~$G$ so that each vertex $t \in H$
is embedded to some $u \in U$.
\end{lemma}

\proof
We may assume without loss of generality that $|U \cap V_1| \geq \lambda m$. If
$t_1 \notin H$, then add $t_1$ to $H$, so now we have $|H| \leq \delta n/6k$.
Moreover, the new $\mathcal{P}_k(H)$ is the union of the old $\mathcal{P}_k(H)$ and
$\{t_1\}$. So now
\begin{equation}\label{newPkH}
|\{x \in T \colon 1 \leq d(x, \mathcal{P}_k(H)) \leq k^3\}| \leq \delta n+\Delta^{k^3+1}\le 3\delta n/2.
\end{equation}
Also, introduce a new constant $\eps'$ with
$\eps \ll \eps' \ll d$. To begin, for each $i$
choose disjoint sets $X_i, Y_i \subseteq V_i$ such that
\begin{itemize}
\item $|X_i| = (1+\alpha/2)m$ and $|Y_i| = 3\lambda m/4 \leq \alpha m/4$,
\item every vertex of $X_i \cup Y_i$ has at least $d\lambda m/2$ inneighbours in
$Y_{i-1}$ and at least $d \lambda m/2$ outneighbours in $Y_{i+1}$, and
\item $Y_1 \subseteq U \cap V_1$.
\end{itemize}
The existence of such sets can be shown by a standard regularity argument.
Indeed, choose disjoint sets $X'_i,Y'_i\subseteq V_i$ such that $|X'_i|=(1+\alpha/2+d^2)m$,
$|Y'_i|=(3\lambda /4+d^2)m$ and $Y'_1\subseteq U\cap V_1$. Then both $G[Y'_{i-1}\to X'_i\cup Y'_i]$
and $G[X'_i\cup Y'_i\to Y'_{i+1}]$ are $2\eps/\lambda$-regular\COMMENT{don't get
$4\eps/3\lambda$-regular since we start with clusters of size $(1+\alpha)m$, also only know that
$\alpha\le 1/2$}
with density at least $3d/4$.
So all but at most $9\eps m/\lambda\le d^2 m$ vertices\COMMENT{get
$4\eps/\lambda(1+\frac{1}{4}+\frac{3\lambda}{4}+2d^2)m\le
(\frac{5\eps}{\lambda}+3\eps+\frac{8\eps d^2}{\lambda})m$} in $X'_i\cup Y'_i$ have
at least $9d\lambda m/16$ inneighbours in $Y'_{i-1}$ and
at least $9d\lambda m/16$ outneighbours in $Y'_{i+1}$. Delete $d^2 m$ vertices from $X'_i$
and $d^2 m$ vertices from $Y'_i$ including these $d^2 m$ vertices of small degree (for each
$i\in [k]$). Then the sets $X_i$ and $Y_i$ thus obtained from $X'_i$ and $Y'_i$ are as
desired.

Each vertex of $\mathcal{P}_k(H)$, and every child of any such vertex, will be
embedded in the sets~$Y_i$, whilst the remaining vertices of $T$ will be
embedded in the sets $X_i$. Observe that by Proposition~\ref{size_of_p(h)},
$|\mathcal{P}_k(H)| \leq 3k|H| \leq \delta n/2$. Moreover, (\ref{newPkH})
implies that there are at most $3\delta n/2$ children of vertices of
$\mathcal{P}_k(H)$ outside $\mathcal{P}_k(H)$. So at most $2 \delta n = d\lambda
m/4$ vertices will be embedded in the sets~$Y_i$.

Next, let $T_1, \dots, T_r$ be the component subtrees of $T[\mathcal{P}_k(H)]$
and $T-\mathcal{P}_k(H)$. So each vertex of $T$ lies in precisely one of the
$T_i$. Let $T^{\textrm{con}}$ be the tree obtained by contracting each $T_i$ to a single
vertex $i$. We may assume the $T_i$ were labelled so that $t_1 \in T_1$ and $1,
2, \dots, r$ is an ancestral order of the vertices of
$T^{\textrm{con}}$. Then let
\begin{align*}
J &=  \{i : T_i \textrm{ is a component subtree of } T[\mathcal{P}_k(H)]\},
\\L &=  \{i : T_i \textrm{ is a component subtree of }T-\mathcal{P}_k(H)\textrm{
and } |T_i| \geq \sqrt{n}\},
\\Q &=  \{i : T_i \textrm{ is a component subtree of }T-\mathcal{P}_k(H)\textrm{
and } |T_i| < \sqrt{n}\}.
\end{align*}
Note that each vertex of $H$ lies in some $T_i$ such that $i \in J$. For each $i
> 1$, $T_i$ contains precisely one vertex with a neighbour in some $T_j$ with $j
< i$. (Furthermore, if $i \in L \cup Q$ then this $j$ must belong to $J$.) Let
$t_i$ be this vertex, then the children of vertices of $\mathcal{P}_k(H)$ which
are not in $\mathcal{P}_k(H)$ are precisely the vertices $t_i$ for $i \in L \cup
Q$. For each~$i$ let $T_i^{\textrm{far}}$ be the set of vertices $x \in T_i$ with $d(t_i,
x) \geq k^3$. Then\COMMENT{Since if $x \in V(T_i) \sm T^{\textrm{far}}$ for some $i \in L
\cup Q$ then $d(x, t_i) < k^3$, and so $x$ is at distance at most $k^3$ from a
member of $\mathcal{P}_k(H)$.}
\begin{equation} \label{eq:not_many_far_vs}
\sum_{i \in L \cup Q} |V(T_i) \sm T^{\textrm{far}}_i| \leq 3\delta n/2
\end{equation}
by~(\ref{newPkH}). Finally, for each $i$ let $T_i^\leq = T[V(T_1) \cup
\dots \cup V(T_{i})]$, so $T_i^\leq$ is the graph formed from the union of
$T_1, \dots, T_i$ by also adding the edges between $T_1, \dots, T_i$.

We shall use a randomised algorithm to embed the vertices of $T$ in $G$. At each
time $\tau$ this algorithm will embed the vertices of $T_\tau$. Indeed, if $\tau
\in J$, we will use Lemma~\ref{embed_p(h)} to embed~$T_\tau$ in the sets $Y_i$
so that the vertices of $H \cap V(T_\tau)$ are embedded in $Y_1 \subseteq U$. If
$\tau \in L$, we will use Lemma~\ref{cluster_cycle_bounded_deg} to embed
$T_\tau$ in the sets $X_i$ (except for the vertex $t_\tau$, which will be
embedded in some $Y_i$) so that approximately equally many vertices of $T_\tau$
are embedded in each set $X_i$. Finally, if $\tau \in Q$ we will use
Lemma~\ref{embed_small_trees} to randomly embed $T_\tau$ in the sets $X_i$
(again with the exception of the vertex $t_\tau$, which will be embedded in some
$Y_i$) so that the expected number of vertices of $T_\tau^{\textrm{far}}$ embedded in each
set $X_i$ is approximately equal. Together the embeddings of each $T_i$ in $G$
will form an embedding of $T$ in $G$ such that every vertex of $H$ is embedded
in $U$, as desired. At any time $\tau$ we will be able to choose the desired
embedding of $T_\tau$ unless there are insufficient vertices remaining
unoccupied in one of the sets $X_i$. We shall show that this is unlikely to
happen for any $i$, and hence that with positive probability the algorithm will
find a copy of $T$ in $G$, proving the lemma.

\medskip \noindent {\bf Tree Embedding Algorithm.}

At time $\tau = 1$, we wish to embed $T_1$ in $G$. Recall that we ensured that
$t_1 \in H$, so $1 \in J$. We shall embed $T_1$ in $Y_1$. Indeed, $|Y_1| = 3\lambda
m/4$, and $|T_1| \leq |\mathcal{P}_k(H)| \leq \delta n/2 = d \lambda m/16$, and so
$Y_1$ contains a copy of $T_1$ by Theorem~\ref{best_so_far}. Choose
(deterministically)\COMMENT{I've put (deterministically) after each of the
deterministic choices to make clear they are non-random. What I mean is that the
algorithm will always make the same choice here, given the same input - this can
be arranged by, for example, declaring (outside the algorithm) a nominated copy
of $T$ in $G$ for any $T$ and $G$ such that $G$ contains a copy of $T$, and then
choosing the nominated copy whenever a choice is made.} such a copy, and embed
each vertex of $T_1$ to the corresponding vertex in this copy.

So after completing the first step, the algorithm will
have obtained an embedding of $T_1 = T_1^\leq$ in $G$ such that any vertex of $H
\cap V(T_1^\leq)$ is embedded in~$Y_1$, and only vertices of $\mathcal{P}_k(H)$
and their children have been embedded in the sets $Y_i$.

At a given time $\tau >1$ we may therefore suppose that the
algorithm has found an embedding of $T_{\tau-1}^\leq$ in $G$ so that each vertex
of $H \cap V(T_{\tau-1}^\leq)$ is embedded in $Y_1$, and only vertices of $\mathcal{P}_k(H)$
and their children have been embedded in the sets $Y_i$. (Recall that this implies that
at most $d\lambda m/4$ vertices are embedded in the sets $Y_i$.) We wish to extend this
embedding to include $T_\tau$, and we do this by the following steps.
\begin{itemize}
\item For each $i$ let $X^\tau_i$ and $Y^\tau_i$ consist of the unoccupied
vertices of $X_i$ and $Y_i$ respectively. If $|X_i^\tau| < |T_\tau|/k+\alpha
m/4$ for some $i$, then terminate the algorithm with failure. So we may assume
that $|X_i^\tau| \geq |T_\tau|/k+\alpha m/4$ for each $i$. Also, since at most
$d\lambda m/4$ vertices have been embedded in the sets $Y_i$, every vertex of
$X_i \cup Y_i$ must have at least $d \lambda m/4$ inneighbours in $Y_{i-1}^\tau$
and at least $d \lambda m/4$ outneighbours in $Y_{i+1}^\tau$.
\item By definition, $t_\tau$ is the unique vertex of $T_\tau$ with a neighbour
which has already been embedded. Let $t'_\tau$ be this neighbour, and let
$v'_\tau$ be the vertex to which $t'_\tau$ was embedded. Also let $V_j$ be the
cluster into which $t_\tau$ should be embedded so that the edge between $t_\tau$
and $t'_\tau$ is embedded canonically. Then $v'_\tau$ has at least $d \lambda
m/4$ neighbours in $Y_j^\tau$, and so by a standard regularity argument, we may
choose some such neighbour $v_\tau \in Y_j^\tau$ which has at least $\alpha
dm/8$ outneighbours in $X_{j+1}^\tau$ and at least $\alpha dm/8$ inneighbours in
$X_{j-1}^\tau$.
\item Now, if $\tau \in L$, for each $i$ consider a set $Z_i^\tau \subseteq X_i^{\tau}$
of size $(1+\alpha/8)|T_\tau|/k$ chosen uniformly at random and independently of all
other choices. We can do this since $(1+\alpha/8)|T_\tau|/k \leq |T_\tau|/k
+\alpha m/8 \leq |X_i^\tau|$ for each $i \in [k]$. Then since $G[X_1^{\tau} \cup
\dots X_k^\tau]$ is a $(16\eps/\alpha)$-regular $d/2$-dense cycle of cluster
tournaments, by Lemma~\ref{regular_restriction} $G[Z_1^\tau, \dots, Z_k^\tau]$
is an $\eps'$-regular $d/4$-dense cycle of cluster tournaments with
probability $1-o(1)$. Also with probability $1-o(1)$, $v_\tau$ has at least
$\alpha
d|T_\tau|/16k$ outneighbours in $Z_{j+1}^\tau$ and at least $\alpha d|T_\tau|/16k$
inneighbours in
$Z_{j-1}^\tau$. So we may choose (deterministically) sets $Z_i^\tau$ satisfying these
two properties. Now delete a single vertex (chosen
arbitrarily) from $Z_j^\tau$, and replace it by $v_\tau$, and let $G^\tau$ be
the restriction of~$G$ to the new $Z_1^\tau, \dots, Z_k^\tau$. Then $G^\tau$ is
a
$(2\eps')$-regular $(d/8)$-dense cycle of cluster tournaments with clusters of
size $(1+\alpha/8)|T_\tau|/k$. So by
Lemma~\ref{cluster_cycle_bounded_deg} $G^\tau$ contains a copy of $T_\tau$ with
at most $(1+\alpha/8) |T_\tau|/k$ vertices of $T_\tau$ embedded in each $X_i$,
and with~$t_\tau$ embedded to $v_\tau$. Embed each vertex of $T_\tau$ to the
corresponding vertex in this copy.
\item If instead $\tau \in Q$, then arbitrarily choose $Z_j^\tau \subseteq X_j^\tau\cup \{v_\tau\}$ of size $\alpha 
m/8$ with $v_\tau \in Z_j^\tau$, and for each $i \neq j$ choose $Z^\tau_i 
\subseteq X_i^\tau$ of size $\alpha m/8$ uniformly at random and independently of 
all other choices. Then $G^\tau := 
G[Z_1^\tau, \dots, Z_k^\tau]$ is a $(16\eps/\alpha)$-regular $d/2$-dense cycle of 
cluster tournaments. Also, with probability $1-o(1)$, $v_\tau$ has at least 
$\alpha^2 dm/128$ outneighbours in $Z_{j+1}^\tau$ and at least $\alpha^2 dm/128$ 
inneighbours in $Z_{j-1}^\tau$, so we may fix (deterministically) our choices of 
the $Z_i^\tau$ such that this event holds. Then by Lemma~\ref{embed_small_trees} 
the set of copies of
$T_\tau$ in $G^\tau$ such that $t_\tau$ is embedded to $v_\tau$ is non-empty,
and furthermore there exists a probability distribution on this set so that if a
copy is chosen according to this distribution, then the expected number of
vertices of $T_\tau^{\textrm{far}}$ embedded in each $Z^\tau_i$ is at most
$(1+\sqrt{\eps})|T_\tau^{\textrm{far}}|/k$. Choose (deterministically) such a
distribution, and choose randomly such a copy according to this distribution.
Embed each vertex of $T_\tau$ to the corresponding vertex in this copy.
\item Finally, if $\tau \in J$, then since $v'_\tau$ has at least $d\lambda m/4$
neighbours in $Y_j^\tau$, we may choose sets $Z_1^\tau \subseteq Y_1^\tau,
\dots, Z_k^\tau \subseteq Y_k^\tau$, each of size $d\lambda m/4$, so that every
vertex of $Z_j^\tau$ is a neighbour of $v'_\tau$. Let $G^\tau$ be the
restriction of $G$ to the sets $Z_i^\tau$; then $G^\tau$ is a $(8\eps/d
\lambda)$-regular $(d/2)$-dense cycle of cluster tournaments. Since $|H| \leq
\delta n/6k = d\lambda m/48k$, by Lemma~\ref{embed_p(h)}, $G^\tau$ contains a
copy of $T_\tau$, with vertex $t_\tau$ embedded in $Y_j^\tau$, and with every
vertex of $H \cap V(T_\tau)$ corresponding to a vertex of $Y_1^\tau$. Embed each
vertex of $T_\tau$ to the corresponding vertex in this copy.
\item
In either case, we have extended the embedding of $T_{\tau-1}^\leq$ in $G$ to an
embedding of $T_{\tau}^\leq$ in $G$, such that every vertex of $H \cap
V(T_\tau^\leq)$ is embedded in $Y_1 \subseteq U$, and only vertices of
$\mathcal{P}_k(H)$
and their children have been embedded in the sets $Y_i$.
\end{itemize}

Since $T_r^\leq = T$, if the algorithm does not terminate with failure then at
time $r$, after embedding $T_r$ it will have obtained an embedding of $T$ in $G$
so that every vertex of $H$ is embedded in $U$, as desired. At this point the
algorithm terminates with success.
\medskip

It remains to show that with positive probability this algorithm will not
terminate with failure before embedding $T_r$. Suppose first that $\sum_{j \in
Q} |T_j| < \alpha m/8$. Then for any $i \in [k]$ and at any time $\tau$, the
number of vertices embedded in $X_i$ is at most
$$\frac{1+\alpha/8}{k} \sum_{\substack{j \in L \\ j < \tau}} |T_j| +
\sum_{\substack{j \in Q \\ j < \tau}} |T_j|
\leq \frac{(1+\alpha/8)(n-|T_\tau|)}{k} + \frac{\alpha m}{8}
< \left(1+\frac{\alpha}{4}\right)m - \frac{|T_\tau|}{k}$$
and so $|X_i^\tau| \geq |T_\tau|/k + \alpha m/4$. Therefore the algorithm cannot
terminate with failure at any point. So we may assume that $\sum_{j \in Q} |T_j|
\geq \alpha m/8$.

Let $OUT$ be the set of all possible courses of the algorithm until termination.
Since the only random choices made by the algorithm are the choices of where to
embed the $T_i$ for each $i \in Q$, any possible course of the algorithm
$\mathcal{C} \in OUT$ can be uniquely described by the embeddings $f_i$ of $T_i$
into $G$ for each $i \in Q$ such that the algorithm does not terminate before
embedding $T_i$. So we may define a probability space with sample space $OUT$
where for any $\mathcal{C} \in OUT$, $\mathbb{P}(\mathcal{C})$ is defined to be
the probability that the algorithm takes course $\mathcal{C}$. So
\begin{align*}
\mathbb{P}(\mathcal{C}) = \prod_{j \in Q} \mathbb{P}(F_j \mid F_i \colon i < j,
i \in Q).
\end{align*}
where $F_j$ denotes the event that $f_j$ is the embedding of $T_j$ into $G$, if
$T_j$ is embedded at some point during $\mathcal{C}$, and is taken to be true
otherwise.

Now, we define the random variable $W_j^i$ in this probability space as follows.
For any $\mathcal{C} \in OUT, j \in Q$ and $i \in [k]$, let
\begin{equation*}
W_j^i(\mathcal{C}) =
\begin{cases}
\frac{\#\textrm{ of vertices from }T_j^{\textrm{far}} \textrm{ embedded in
}X_i}{\sqrt{n}} & \text{if $T_j$ is embedded during $\mathcal{C}$,}
\\ \frac{|T_j^{\textrm{far}}|}{k\sqrt{n}} & \text{otherwise.}
\end{cases}
\end{equation*}
Since $|T_j^{\textrm{far}}| \leq |T_j| < \sqrt{n}$ for each $j \in Q$,
$W_j^i$ is a well-defined function from $OUT$ to $[0,1]$, and so is a
well-defined random variable in our probability space.

For any $j \in Q$ and $\mathcal{C}_a, \mathcal{C}_b \in OUT$, let
$\mathcal{C}_a \sim_j \mathcal{C}_b$ if and only if $\mathcal{C}_a$ and
$\mathcal{C}_b$ share the same course before time~$j$ (i.e. they embed $T_1,
\dots, T_{j-1}$ identically) or $T_j$ is not embedded at any point in either
$\mathcal{C}_a$ or~$\mathcal{C}_b$. Then $\sim_j$ is an equivalence relation on
$OUT$ (since if two courses agree up to time $j-1$, then at time $j$ either they both terminate with failure or they both successfully embed $T_j$). For any equivalence class $\mathcal{C}^*$ of $\sim_j$ other than the
class of $\mathcal{C}$ for which $T_j$ is not embedded, every $\mathcal{C} \in
\mathcal{C}^*$ shares the same course before time~$j$. So for each $\mathcal{C}
\in \mathcal{C}^*$, the same probability distribution on the set of copies of
$T_j$ will have been chosen at time~$j$, and a copy will then have been chosen
according to this distribution. So further partition $\mathcal{C}^*$ into
$\mathcal{C}^*_1, \dots, \mathcal{C}^*_a$ by this choice, so courses
$\mathcal{C}, \mathcal{C}' \in \mathcal{C}^*$ are in the same $\mathcal{C}^*_s$
if and only if $T_j$ is embedded identically in $\mathcal{C}$ and
$\mathcal{C}'$.
Now
$$\mathbb{E}(W_j^i \mid \mathcal{C}^*) = \sum_s \mathbb{E}(W_j^i \mid
\mathcal{C}^*_s)\mathbb{P}(\mathcal{C}^*_s \mid \mathcal{C}^*),$$
but every member of $\mathcal{C}^*_s$ embeds $T_j$ identically, so
$\mathbb{E}(W_j^i \mid \mathcal{C}^*_s)$ is simply the number of vertices of
$T_j^{\textrm{far}}$ embedded in $X_i$ in this common embedding, divided by $\sqrt{n}$.
Also, $\mathbb{P}(\mathcal{C}^*_s \mid  \mathcal{C}^*)$
is the probability that this embedding of $T_j$ is chosen when the random choice
of the embedding of~$T_j$ is made. So by our (deterministic) choice of the
probability distribution on the copies of $T_j$ in $G$,
\begin{equation} \label{eq:expected_z}
\mathbb{E}(W_j^i \mid \mathcal{C}^*) \leq \frac{(1 +
\sqrt{\eps})|T_j^{\textrm{far}}|}{k\sqrt{n}}.
\end{equation}
If instead $\mathcal{C}^*$ is the class of all $\mathcal{C}$ such that $T_j$ is
not embedded in $\mathcal{C}$, then $\mathbb{E}(W_j^i \mid \mathcal{C}^*) =
|T_j^{\textrm{far}}|/k\sqrt{n}$ by definition, and so (\ref{eq:expected_z}) holds in this
case also.

Now, for any equivalence class $\mathcal{C}^*$ other than the class in which
$T_j$ is not embedded, the embeddings of $T_1, \dots, T_{j-1}$ are identical
amongst the members of $\mathcal{C}^*$, and so
$$\mathbb{E}(W_j^i\mid \mathcal{C}^*, W_s^i \colon s \in Q, s < j) =
\mathbb{E}(W_j^i\mid \mathcal{C}^*).$$
Clearly this equality also holds for the class $\mathcal{C}^*$ in which $T_j$ is
not embedded, and so for any $i \in [k]$,
\begin{align*}
\mathbb{E}(W_j^i \mid W_s^i \colon s \in Q, s < j) &\leq \max_\mathcal{C^*}
\mathbb{E}(W_j^i \mid \mathcal{C}^*, W_s^i \colon s \in Q, s < j) \leq \frac{(1
+ \sqrt{\eps})|T_j^{\textrm{far}}|}{k\sqrt{n}}.
\end{align*}
Since $\sum_{j \in Q}|T_j| \geq \alpha m/8$, by Lemma~\ref{azuma}, for any $i$
the probability that \begin{equation} \label{eq:sum_z_small}
\sum_{j\in Q} W_j^i \leq \frac{(1+\alpha/8)\sum_{j \in Q} |T_j|}{k\sqrt{n}}
\end{equation}
does not hold decreases exponentially with $n$. So with probability $1-o(1)$,
(\ref{eq:sum_z_small}) holds for each $i \in [k]$.

To finish the proof, we show that if (\ref{eq:sum_z_small}) holds for each $i
\in [k]$, then the algorithm cannot terminate with failure, and will therefore
successfully embed $T$ in $G$ as desired. Indeed, the algorithm will only
terminate with failure if at some time $\tau$ we have $|X_i^\tau| <
|T_\tau|/k+\alpha m/4$ for some $i$. But for any $i \in [k]$ and any time
$\tau$, only vertices from subtrees $T_s$ such that $s \in L \cup Q$ and $s <
\tau$ have been embedded in $X_i$ before time $\tau$. So the number of vertices
embedded in~$X_i$ before time $\tau$ is at most\COMMENT{For the second
line, if $\tau \in L$ then we have taken off $|T_\tau|/k$ and have added on more
than this (in the first term). Alternatively, if $\tau \in Q \cup J$ then we
have taken off $|T_\tau|/k \leq \delta n/2$ and have added on $\delta n/2$. We also
use $2\delta n \leq \alpha m/8$.}
\begin{align*}
\frac{(1+\alpha/8)}{k} \sum_{s \in L \sm \{\tau\}}|T_s| &+ \sum_{s \in Q \sm
\{\tau\}} |V(T_s) \sm T_s^{\textrm{far}}| + \sqrt{n}\sum_{s \in Q \sm \{\tau\}} W_s^i
\\ &\stackrel{(\ref{eq:not_many_far_vs})}{\leq} \frac{(1+\alpha/8)}{k} \sum_{s
\in L} |T_s| + \frac{3\delta n}{2} + \sqrt{n} \sum_{s \in Q} W_s^i + \frac{\delta n}{2} -
\frac{|T_\tau|}{k}
\\ &\stackrel{(\ref{eq:sum_z_small})}{\leq} \frac{(1+\alpha/8)}{k} \sum_{s \in
L \cup Q} |T_s| + 2\delta n  - \frac{|T_\tau|}{k}
\leq (1+\frac{\alpha}{4})m - \frac{|T_\tau|}{k}.
\end{align*}
To see that the second line holds, note that $|T_\tau|/k < \sqrt{n}/k<\delta n/2$ whenever
$\tau \in Q$ and $|T_\tau|/k < |\mathcal{P}_k(H)|<\delta n/2$ whenever
$\tau \in J$.
So if (\ref{eq:sum_z_small}) holds, then at any time $\tau$ and for any $i \in
[k]$, $|X_i^\tau| \geq |T_\tau|/k+\alpha m/4$, and so the algorithm succeeds.
This completes the proof of Lemma~\ref{embed_with_restrictions}.
\endproof

\subsection{Proof of Lemma~\ref{robust_expander_case}.} \label{subsec:rob_exp_proof2}
We can now give the proof of Lemma~\ref{robust_expander_case}, which will proceed 
as follows. We shall apply Lemma~\ref{extendedtree} to find a subtree $T_{\textrm{ext}}$ of $T$
and a subset $H \subseteq V(T_{\textrm{ext}})$. Then we shall find a cluster cycle 
$\mathcal{C}$ in~$G$ such that $|\mathcal{C}|$ is slightly larger than 
$|T_{\textrm{ext}}|$. We then embed $T_{\textrm{ext}}$ into $\mathcal{C}$ using 
Lemma~\ref{embed_with_restrictions}, restricting $H$ to a set $U$ of vertices 
of~$\mathcal{C}$ which have many inneighbours and outneighbours outside 
$\mathcal{C}$. Finally we shall use this property of $U$ to embed the vertices of 
$T - T_{\textrm{ext}}$ in $V(G) \sm V(\mathcal{C})$
and thereby complete the embedding.

If $\alpha \geq 1/2$, then $G$ contains a copy of $T$ by Theorem~\ref{best_so_far}.
So we may assume that $\alpha < 1/2$.
We begin by introducing new
constants $\Delta^*, M, M', \delta, \eps, d$ and $\Delta$ with
\begin{align*}
1/n \ll 1/\Delta^* \ll& 1/M \ll 1/M' \ll \eps  \ll d \ll \mu \ll \nu \ll \eta \ll 1/\Delta \ll \alpha.
\end{align*}
Then Lemma~\ref{KOT_combined} implies that $G$ contains an $\eps$-regular 
$d$-dense cycle of cluster tournaments on clusters $V_1, \dots, V_k$, where 
$M' \leq k \leq M$ and each cluster has equal size between $(2+\alpha)m$ and 
$(2+2\alpha)m$, where $m = n/k$. Also let
$$
\delta := d\alpha^2/16000k.
$$
Remove vertices from each $V_i$ to obtain a 
$2\eps$-regular $d/2$-dense cycle of cluster tournaments $G'$ on clusters $V_1', \dots, V_k'$ each of size $(2+\alpha) m$.

Let $T_c$ be the core tree of $T$ with parameter $\Delta$, and choose any 
vertex~$t_1 \in T_c$ as the root of $T$. Then by Lemma~\ref{extendedtree} (applied with $\omega = \delta \alpha /50$), we may 
choose a subtree $T_{\textrm{ext}}$ of $T$ and a subset $H \subseteq V(T_{\textrm{ext}})$ satisfying 
the following properties.
\begin{itemize}
\item[(i)] $T_c \subseteq T_{\textrm{ext}}$.  
\item[(ii)]  $\Delta(T_{\textrm{ext}}) \leq \Delta^*$.
\item[(iii)]  For any edge $e$ between $V(T - T_{\textrm{ext}})$ and $V(T_{\textrm{ext}})$, the endvertex of $e$ in $V(T_{\textrm{ext}})$ lies in~$H$.
\item[(iv)]  The number of
vertices $v \in T_{\textrm{ext}}$ which satisfy $1 \leq d(v, \mathcal{P}_k(H)) \leq
k^3$ is at most $\delta \alpha n/50$.
\item[(v)]  $|H| \leq n/\Delta^{k^{50/\delta \alpha}} \leq \delta \alpha n/350k$.
\end{itemize}

Let $T_1^+, \dots, T_r^+$ and $T_1^-, \dots, T_s^-$ be the component
subtrees of $T - T_{\textrm{ext}}$. Each $T_i^+$ and $T_i^-$ will contain precisely one
vertex, $v_i^+$ or $v_i^-$ respectively, with a neighbour in $T_{\textrm{ext}}$. Label the
$T_1^+, \dots, T_r^+$ and $T_1^-, \dots, T_s^-$ so that each $T_i^+$ contains
$v_i^+$ with an inneighbour in $T_{\textrm{ext}}$, and each $T_i^-$ contains $v_i^-$ with
an outneighbour in $T_{\textrm{ext}}$. By (i) and Proposition~\ref{core_tree_props}(iv) each
$T_i^+$ and each $T_i^-$ contains at most $n/ \Delta$ vertices. Let $x = 
|T_{\textrm{ext}}|$, let $y = |T_1^+ \cup \dots \cup T_r^+|$ and let $z = |T_1^- \cup \dots 
\cup T_s^-|$, so $x+y+z=n$. 

Then all but at most
$2y+x+\alpha n/2$ vertices of $G$ have at least $y+x/2+\alpha n/4$
outneighbours in $G$, and all but at most $2z+x+\alpha n/2$ vertices of $G$ have
at least $z + x/2 + \alpha n/4$ inneighbours in $G$. So at least $2(1+\alpha)n
- 2y
- 2z - 2x - \alpha n = \alpha n$ vertices of $G$ satisfy both of these
conditions. Choose any $\alpha n/8$ of these vertices to form $U_0$. Then $|U_0|
= \alpha n/8$, and each $v \in U_0$ has
at least $y+ x/2 + \alpha n/8$ outneighbours outside $U_0$ and at least $z +
x/2+ \alpha n/8$ inneighbours outside $U_0$.

Suppose first that $x \geq \alpha n/50$.  From each
cluster $V'_i$ of $G'$ choose a set $X_i$ of $x(1+\alpha /2)/k$ vertices
uniformly at random, and let $X = X_1 \cup \dots \cup X_k$. Then $|X| = 
x(1+\alpha/2)$, and for any  
single vertex $u \in G'$, the probability that~$u$ is included in
$X$ is equal to $x/2n$. So by Proposition~\ref{chernoff}, with probability
$1-o(1)$ the set $U := X \cap U_0$ satisfies $|U| \geq \alpha x/20 \geq
\alpha^2n/1000$. Also, for any vertex $v \in U$, the expected number of 
outneighbours
of~$v$ outside $X$ is at least
\begin{align*}
\left(y + \frac{x}{2} + \frac{\alpha n}{8}\right) \left(1 - \frac{x}{2n}\right)
&=  y + \frac{x}{2} - \frac{xy}{2n} - \frac{x^2}{4n} + \left(1-
\frac{x}{2n}\right)
\frac{\alpha n}{8}
\\ &\geq y + x\left(\frac{1}{2} - \frac{y+x}{2n}\right) + \frac{\alpha n}{16}
\geq y + \frac{\alpha n}{16}.
\end{align*}
A similar calculation shows that for each $v \in U$, the expected number of
inneighbours of $v$ outside $X$ is at least $z + \alpha n/16$. So by
Proposition~\ref{chernoff} we find that with probability $1-o(1)$, every vertex
$v \in U$ has at least $y + \alpha n/20$ outneighbours outside $X$ and at least
$z + \alpha n/20$ inneighbours outside $X$. Fix a choice of $X$ such that both
these events of probability $1-o(1)$ occur. 

Since every vertex of $U$ has either at least\COMMENT{
$|X| = x(1+\alpha/2)$ so $(2(1+\alpha)n - |X|)/2 = n + \alpha n - x/2 -
\alpha x/4 = x/2 + y + z + \alpha n - \alpha x/4$}
$(2(1+\alpha)n - |X|)/2 \geq y + z + x/2 +\alpha n/2 \geq y+z+\alpha n/2$
inneighbours outside $X$ or at
least $y+z+\alpha n/2$
outneighbours outside $X$, we may choose a set $U' \subseteq U$ of size $|U'|
\geq
|U|/2 \geq \alpha^2 n/2000$ such that either
\begin{itemize}
\item[(a)] every $v \in U'$ has at least $y + \alpha n/20$ outneighbours outside
$X$ and at least $y+z+\alpha n/20$ inneighbours outside $X$, or
\item[(b)] every $v \in U'$ has at least $y+z+\alpha n/20$ outneighbours
outside $X$ and at least $z + \alpha n/20$ inneighbours outside $X$.
\end{itemize}

So $G'[X]$ is a $(150\eps/\alpha)$-regular $(d/2)$-dense cycle of cluster 
tournaments on clusters $X_1, \dots, X_k$ each of size $(1+\alpha/2)x/k$, and $U' 
\subseteq X_1 \cup \dots X_k$ has size $|U'| \geq \alpha^2 n/2000 \geq \alpha^2 x/2000$. Also by 
(ii), (iv) and (v) we know that $T_{\textrm{ext}}$ is a directed tree on $x$ vertices 
rooted at $t_1$ and with $\Delta(T_{\textrm{ext}}) \leq \Delta^*$, and also that $H 
\subseteq V(T_{\textrm{ext}})$ satisfies $|H| \leq \delta \alpha n/350k \leq \delta x/7k$ 
and $|\{t \in T_{\textrm{ext}}: 1 \leq d(t, \mathcal{P}_k(H)) \leq k^3\}| \leq \delta 
\alpha n/50 \leq \delta x$. So by Lemma~\ref{embed_with_restrictions} (with 
$\Delta^*$ in place of $\Delta$ and $\alpha^2/2000$ in place of $\lambda$), we may embed $T_{\textrm{ext}}$ in $G'[X]$ so that every 
vertex of $H$ is embedded to a vertex of $U'$.

Now suppose instead that $x < \alpha n/50$. Then, since every vertex $v$ of $G$
has either $d^+(v) \geq (1+\alpha)n-1 \geq y+z+\alpha n$ or $d^-(v) \geq (1+\alpha)n-1 \geq y+z+\alpha n$, we can
choose a set $U' \subseteq U_0$ of size $|U'| \geq \alpha n/16$ which satisfies
either (a) or (b) (with $X := U'$). Since $|T_{\textrm{ext}}| = x
< \alpha n/50 \leq |U'|/3$, and $G[U']$ is a tournament, by
Theorem~\ref{best_so_far} we may embed $T_{\textrm{ext}}$ in $G[U']$, so in particular
every vertex of $H$ is embedded to a vertex of $U'$.

In either case, let $V_{\textrm{ext}}$ be the set of vertices of $G$ to which $T_{\textrm{ext}}$
is embedded. We may now complete the embedding of $T$ in $G$. If $U'$ satisfies
(a), then we first proceed through the trees~$T^+_i$ in turn. For each $T^+_i$,
let $u^+_i$ be the inneighbour of $v^+_i$ in $T_{\textrm{ext}}$ (so $u^+_i \in H$ by (iii)). Then
$u^+_i$ has been
embedded to some vertex $v \in U'$. This $v \in U'$ has at least $y + \alpha
n/20$ outneighbours outside $V_{\textrm{ext}}$, of which at most $y$ have been used for
embedding the trees $T^+_j$ for $j<i$. So there are at least $\alpha n/20$
outneighbours
of $v$ outside $V_{\textrm{ext}}$ available to embed $T^+_i$, and so since $|T^+_i| \leq
n/\Delta \leq \alpha n/60$, by Theorem~\ref{best_so_far} we can embed $T^+_i$
among these vertices. In this way we may embed each of the $T^+_i$. We then
proceed through the $T^-_i$ similarly. For each $T^-_i$ let $u^-_i$ be the
inneighbour of $v^-_i$ in $T_{\textrm{ext}}$ (so $u^-_i \in H$ by (iii)). Then $u^-_i$ has been
embedded to some
vertex $v \in U'$. This $v \in U'$ has at least $y + z + \alpha n/20$
inneighbours outside $V_{\textrm{ext}}$, of which at most $y+z$ have been used for
embedding
the trees $T^+_1, \dots, T^+_r$ and the trees $T^-_j$ for $j<i$. So there are at
least
$\alpha n/20$ inneighbours of $v$ outside $V_{\textrm{ext}}$ available to embed $T^-_i$, and
so since $|T^-_i| \leq n/\Delta \leq \alpha n/60$, again by
Theorem~\ref{best_so_far} we can embed $T^-_i$ among these vertices. If $U'$
satisfies (b) we can
embed $T$ similarly, first embedding the $T^-_i$, and then the $T^+_i$. Either
way we have completed the embedding of $T$ in $G$.
\endproof

\section{Embedding trees in an almost-transitive
tournament.}\label{sec:transitive}

A \emph{transitive tournament} is a tournament in which the vertices can be
given a total order so that every edge is directed towards the endvertex which
is greater in this order. It is easy to show that any transitive tournament $G$
on $n$ vertices contains any directed tree $T$ on $n$ vertices, by first showing
that the
vertices of $T$ can be given\COMMENT{
Induction on $n$: the case $n=1$ is trivial, and for $n \geq 1$ choose any edge
$e = t_1 \rightarrow t_2$ of $T$ and delete it from $T$, use the inductive
hypothesis to choose such orders for the vertices of each component of $T$ with
this edge deleted; then let the order of the vertices of $T$ be the vertices of
the component containing $t_1$ (in their order) followed by the vertices of the
component containing $t_2$ (in their order).
} a total order so that every edge is directed towards the endvertex which is
greater in this order, and then embedding each vertex of $T$ to the vertex of $G$
in the corresponding position (in the order of vertices of $G$).

In this section, we shall prove an approximate version of this result, namely
that if a tournament on $(1+\alpha)n$ vertices (for some small $\alpha$) is
sufficiently close to being transitive, then it contains any directed tree on $n$
vertices. To state this lemma precisely, we say that a tournament~$G$ on $n$
vertices is \emph{$\eps$-almost-transitive} if the vertices of $G$ can be given an
order $v_1, \dots, v_n$ so that at most $\eps n^2$ edges are directed against
the ordering of the vertices, that is, they are directed from $v_i$ to $v_j$
where $i > j$.

The proof of this lemma is by a similar method to the proof of
Theorem~\ref{main} in the next section. The approach is that if the lemma is false,
then
there is some $\alpha > 0$ for which the lemma does not hold, and so the infimum
$a_{\textrm{inf}}$ of all $\alpha$ for which the lemma does hold is greater than zero. We
then choose $\alpha$ slightly less than $a_{\textrm{inf}}$ and apply (to a smaller
subtree) the fact that the lemma
holds for any $\alpha' > a_{\textrm{inf}}$ to show that the lemma holds for $\alpha$,
giving a contradiction.

\begin{lemma} \label{almost_trans}
For all $\alpha >0$ there exists $\eps_0>0$ and $n_0 \in \mathbb{N}$ such that
for any $\eps \leq \eps_0$ and any $n \geq n_0$, any $\eps$-almost-transitive
tournament $G$ on at least $(1+\alpha)n$ vertices contains any directed tree $T$
on $n$ vertices.
\end{lemma}

\proof
We consider the set $A$ of all positive values of $\alpha$ such that the lemma
holds. More precisely, $A$ is the set of all positive values of $\alpha$ such
that there exist $\eps_0 >0$ and $n_0 \in \mathbb{N}$ so that for any $n \geq
n_0$ and $\eps \leq \eps_0$, any $\eps$-almost-transitive tournament $G$ on at
least
$(1+\alpha)n$ vertices contains a copy of any directed tree $T$ on $n$ vertices.
So if $\alpha' \in A$ and $\alpha'' > \alpha'$ then $\alpha'' \in A$. Also $2 \in A$
by Theorem~\ref{best_so_far},
and so we may define $a_{\textrm{inf}} = \inf A$, with $0 \leq a_{\textrm{inf}} \leq 2$.
Then for any $\alpha' > a_{\textrm{inf}}$, $\alpha' \in A$. With this definition the
lemma is
equivalent to the statement that $a_{\textrm{inf}} = 0$, so suppose for a contradiction
that $a_{\textrm{inf}} >0$. Let
$$\gamma \ll 1/\Delta \ll a_{\textrm{inf}} \textrm{\hspace{1cm} and \hspace{1cm}}\alpha =
a_{\textrm{inf}} - \gamma,$$
so we may assume that $1/\Delta \ll \alpha$. Then $\alpha+2\gamma > a_{\textrm{inf}}$, so
$\alpha +2\gamma \in A$, and so by definition of $A$ there exist $\eps_0'>0$ and
$n'_0 \in \mathbb{N}$ such that for any $\eps' \leq \eps_0'$ and $n' \geq n'_0$,
any $\eps'$-almost-transitive tournament $G$ on at least $(1+\alpha+2\gamma)n'$
vertices
contains a copy of any directed tree $T$ on $n'$ vertices. Moreover, we may assume that
$\eps_0' \ll \gamma$. Fix such an $\eps'_0$ and
$n_0'$, and let $1/n_0 \ll 1/n'_0, \gamma$ and $\eps_0 \ll \eps'_0$. We
will show that for any $n \geq n_0$ and $\eps \leq \eps_0$, any
$\eps$-almost-transitive tournament $G$ on at least $(1+\alpha)n$ vertices
contains a
copy of any directed tree $T$ on $n$ vertices. It then follows that $\alpha \in A$,
yielding a contradiction and therefore proving the lemma.

So let $\eps \leq \eps_0$ and $n \geq n_0$, let $G$ be an
$\eps$-almost-transitive tournament on at least $(1+\alpha)n$ vertices and let
$T$ be a directed tree on $n$ vertices. If $|G| \geq 3n$, then $G$ contains a a
copy of~$T$ by Theorem~\ref{best_so_far}, and so we may assume that $|G| < 3n$.
Since $G$ is $\eps$-almost-transitive, we may order the
vertices of $G$ as $v_1, \dots, v_{|G|}$ so that at most $\eps |G|^2 \leq 9 \eps
n^2$ edges are directed from $v_j$ to $v_i$ where $i < j$. Now, at most\COMMENT{
Since each such edge is incident to two vertices, the sum over every vertex
$v$ of the number of such edges $v$ is incident to is at most $18 \eps n^2$.}
 $18 \sqrt{\eps} n$ vertices of $G$ are incident to more than $\sqrt{\eps} n$
such edges; let $G'$ be the subgraph of $G$ obtained by deleting these vertices
from $G$, and let $v'_1, v'_2, \dots, v'_{|G'|}$ be the vertices of $G'$ in the
inherited order. Then $G'$ is a tournament on at least $(1+\alpha - 18
\sqrt{\eps})n$ vertices such that for any vertex $v'_i$ there are at most
$\sqrt{\eps} n$ vertices~$v'_j$ for which the edge between $v'_i$ and $v'_j$ is
directed
towards $v'_{\min\{i, j\}}$.

Next, let $T_c$ be the core tree of $T$ with parameter $\Delta$, as defined
in Section~\ref{subsec:core_tree}. We consider three possibilities for $T_c$, in
each case showing that
$T$ can be embedded in $G'$.

\emph{Case 1: Some vertex $t \in T_c$ has $d^+_{T_c}(t) \geq 2$.} Then let $F^-$
be the (possibly empty) forest
consisting of each component subtree $T'$ of $T - t$ such that the edge between
$T'$ and $t$
is directed towards~$t$. Similarly let the component subtrees $T''$ of $T - t$
such that
the edge between $T''$ and $t$ is directed away from $t$ be partitioned into two
forests, $F_1^+$ and $F_2^+$. Since $d^+_{T_c}(t) \geq 2$, by
Proposition~\ref{core_tree_props}(ii) at least two such component subtrees each
contain at least $n/\Delta$ vertices, and so we
may choose
$F_1^+$ and $F_2^+$ so that $|F_1^+|, |F_2^+| \geq n/\Delta$. Note that $|F^-| =
w^-(t)$, and $|F^+_1|+|F^+_2|=w^+(t)$, so in particular $w^+(t) \geq 2n/\Delta$,
and also recall that $w^+(t) +w^-(t) =n-1$.

We first determine where to embed the vertex $t$. For this, let
\begin{equation*}
p :=
\begin{cases} 3\gamma n +\sqrt{\eps}n+1 & \text{if $w^-(t) < \gamma n$,}
\\
(1+\alpha+2\gamma)w^-(t)+\sqrt{\eps}n+1 &\text{if $w^-(t) \geq \gamma n$.}
\end{cases}
\end{equation*}
and embed $t$ to the vertex $v'_p$ of $G'$. This can be done, as we shall see
later that $p < |G'|$. We will embed $F^-$ in the vertices preceding $v'_p$ and
$F_1^+, F_2^+$ in the vertices succeeding $v'_p$ in the vertex ordering of $G'$.
Embedding $F^-$ will be possible because $p$ is a little larger than one might
expect, whereas embedding $F_1^+$ and $F_2^+$ can be done successively, which
will give us enough room for both. Let $S^- = N^-(v'_p) \cap \{v'_1,
\dots, v'_{p-1}\}$, and $S^+ = N^+(v'_p) \cap \{v'_{p+1}, \dots, v'_{|G'|}\}$.
Then~$S^-$ and $S^+$ are disjoint, $|S^-| \geq p - \sqrt{\eps}n - 1$
and\COMMENT{Since $v'_p$ is incident to at most $\sqrt{\eps} n$ edges which go
`the wrong way'.}
$|S^+| \geq |G'| - p - \sqrt{\eps}n$. Next we shall embed $F^-$ in $G'[S^-]$.
Indeed, if $w^-(t) < \gamma n$ then $|S^-| \geq 3 \gamma n$, and so by
Theorem~\ref{best_so_far} we can embed $F^-$ in $G'[S^-]$.
Alternatively, if $w^-(t) \geq \gamma n$, let $n' = w^-(t) \geq n'_0$ and $\eps' =
|G|^2\eps/(n')^2\leq \eps'_0$, then $F^-$ is a forest on $n'$ vertices, and
$G'[S^-]$ is an $\eps'$-almost-transitive\COMMENT{Since at most $\eps|G|^2$
edges are directed `the wrong way' in $G$, and so at most this many are directed
`the wrong way' in $G[S^-]$.}
 tournament on at least $(1+\alpha+2\gamma) n'$ vertices. So by the choice of
$\eps'_0$ and $n'_0$ we can embed~$F^-$ in $G'[S^-]$.

Finally we shall complete the embedding of $T$ in $G'$ by embedding $F^+_1$ and
$F^+_2$ in $G'[S^+]$. Now,
\begin{align*}
|S^+| &\geq  |G'| - p - \sqrt{\eps}n
\\ &\geq (1+\alpha -18\sqrt{\eps})n - (3 \gamma n +
(1+\alpha+2\gamma)w^-(t)+\sqrt{\eps}n +1) - \sqrt{\eps}n
\\ &\geq (1+\alpha)w^+(t) - 5\gamma n - 20\sqrt{\eps}n \geq (1+\alpha)w^+(t) -
6\gamma n.
\end{align*}
Let $n' = |F_1^+|$, so $n'_0 \leq
n/\Delta \leq n'$ and $n' \leq w^+(t)-n/\Delta$, and again let $\eps' = |G|^2\eps/(n')^2$,
so $\eps' \leq \eps'_0$. Then $G'[S^+]$ is an $\eps'$-almost-transitive
tournament on
$|S^+| \geq (1+\alpha)(n'+n/\Delta) - 6\gamma n \geq
(1+\alpha+1/\Delta)n'+(\alpha/\Delta - 6\gamma)n \geq (1+\alpha+2\gamma)n'$
vertices, and so by our choice of $n'_0$ and $\eps'_0$,
we may embed $F_1^+$ in $G'[S^+]$.

Now, let $S^+_{\textrm{rem}}$ consist of the vertices of $S^+$ not occupied by the
vertices of $F_1^+$.
We shall embed $F_2^+$ in $S^+_{\textrm{rem}}$ in a similar manner. Indeed, we now let
$n' = |F_2^+|$, so again $n'_0 \leq n/\Delta \leq n'$, and
again let $\eps' = |G|^2\eps/(n')^2 \leq \eps'_0$. Then\COMMENT{Since
$\alpha/\Delta \geq 8\gamma$ and $n \geq n'$, using $w^+(t) =
|F_1^+|+|F_2^+|$ and $|F_1^+| \geq n/\Delta$,
$|S^+_{\textrm{rem}}| \geq (1+\alpha)w^+(t) - 6\gamma n - (w^+(t)
-|F_2^+|) = |F_2^+|+\alpha w^+(t) - 6\gamma n = (1+\alpha) n'
+ \alpha |F_1^+| - 6\gamma n \geq (1+\alpha) n' + (\alpha/\Delta - 6\gamma)n'$.
}
\begin{align*}|S^+_{\textrm{rem}}| = |S^+| - |F_1^+| &\geq (1+\alpha)w^+(t) - 6\gamma n -
(w^+(t)
-|F_2^+|) \\&= (1+\alpha) n'
+ \alpha |F_1^+| - 6\gamma n \geq (1+\alpha+2\gamma)n',
\end{align*}
 so
$G'[S^+_{\textrm{rem}}]$ is an $\eps'$-almost-transitive tournament on at least
$(1+\alpha+2\gamma)n'$ vertices, and so by our choice of $n'_0$ and $\eps'_0$,
we may embed $F_2^+$ in $G'[S^+_{\textrm{rem}}]$.

\emph{Case 2: Some vertex $t \in T_c$ has $d^-_{T_c}(t) \geq 2$.} Then we may
embed
$T$ in $G'$ by the same method as in Case 1, the main
difference
being that the roles of outdegrees and outneighbours are switched with those of
indegrees and inneighbours.

\emph{Case 3: $T_c$ is a directed path (possibly consisting of just a single
vertex).} Then let $w^+ = w^+(T_c)$ and $w^- = w^-(T_c)$ be as defined in
Section~\ref{subsec:defs}, and partition the vertices of $G'$ into three sets
$S^- = \{v'_{1}, \dots, v'_{w^-+\alpha n/3}\}$, $S = \{v'_{w^-+\alpha n/3+1},
\dots, v'_{|G'| - w^+ - \alpha n/3}\}$ and $S^+ = \{v'_{|G'| -
w^+ - \alpha n/3+1}, \dots, v'_{|G'|}\}$.
Then since $w^+ +w^- + |T_c| = n$, we know that
$|S| = |G'| - w^- - w^+ - 2\alpha n/3 \geq |T_c|$. Therefore by
Theorem~\ref{dir_path_embed} we may embed $T_c$ in $G'[S]$. Now, let $T_1^+,
\dots,
T_r^+$ be the component subtrees of $T - T_c$ such that the edge between $T_i^+$
and $T_c$ is directed towards $T_i^+$, and for each $i$ let $t_i^+ \in T_c$ be
the vertex of $T_c$ to which this edge is incident, and let $v_i^+$ be the
vertex of $G'$ to which $t_i^+$ was embedded. Similarly, let $T_1^-, \dots,
T_s^-$ be the component subtrees of $T - T_c$ such that the edge between $T_i^-$
and $T_c$ is directed towards $T_c$, let $t_i^-$ be the vertex of $T_c$ to which
this edge is incident, and let $v_i^-$ be the vertex of $G'$ to which $t_i^-$
was embedded. Then every vertex of $T$ lies in $T_c$ or one of the $T_i^+$ or
$T_i^-$. Furthermore $|T_i^+|, |T_j^-| \leq
n/\Delta$ for each $i$ and $j$ by Proposition~\ref{core_tree_props}(iv).

We shall complete the embedding of $T$ in $G'$ by greedily embedding each
$T_i^+$ in $N^+(v_i^+) \cap S^+$, and each $T_i^-$ in $N^-(v_i^-) \cap S^-$.
Indeed, suppose we have already embedded $T_1^+, \dots, T_{i-1}^+$, and we now
wish to embed $T_i^+$. Then
$$|N^+(v_i^+) \cap S^+| \geq |S^+| - \sqrt{\eps}n \geq w^+ + \alpha n/3 -
\sqrt{\eps}n \geq w^+ + \alpha n/4.$$
At most $w^+$ of these vertices have
already been occupied by vertices of $T_1^+, \dots, T_{i-1}^+$, and so there
remain at least $\alpha n/4$ available vertices in which to embed $T_i^+$.
Since
$|T_i^+| \leq n/\Delta \leq \alpha n/12$, we may embed $T_i^+$ in these
available vertices by Theorem~\ref{best_so_far}. Continuing in this way we may
embed all of the $T_i^+$, and the $T_i^-$ may be embedded similarly\COMMENT{
Suppose we have embedded $T_1^-, \dots, T_{i-1}^-$ in $N^-(v_i^-) \cap S^-$, and
we now wish to embed $T_i^-$. Then $|N^-(v_i^-) \cap S^-| \geq |S^-| -
\sqrt{\eps}n \geq w^- + \alpha n/3 - \sqrt{\eps}n \geq w^- + \alpha n/4$. At
most $w^-$ of these vertices have already been occupied by vertices of $T_1^-,
\dots, T_{i-1}^-$ (and none have been occupied by any $T_j^+$, as $S^+$ and
$S^-$ are disjoint), and so there remain at least $\alpha n/4$ available
vertices in which to embed $T_i^-$. Since $|T_i^-| \leq n/\Delta \leq \alpha
n/12$, we may embed $T_i^-$ in these available vertices by
Theorem~\ref{best_so_far}.
}, to give us a copy of~$T$ in $G'$.

Any tree in which every vertex has at most one outneighbour and at most one
inneighbour is a directed path. So $T_c$ must fall into at least one of the
three cases, and
so we can find a copy of $T$ in $G'$, and hence in $G$, contradicting our
assumption that $a_{\textrm{inf}} > 0$. So we must have $a_{\textrm{inf}} = 0$, and so the lemma
holds.
\endproof

\section{Proof of Theorem~\ref{main}} \label{proof}

Recall the statement of Theorem~\ref{main}. \medskip

\noindent\textbf{Theorem~\ref{main}}
\emph{Let $\alpha >0$. Then the following properties hold.
\begin{enumerate}
\item There exists $n_0 \in \mathbb{N}$ such that for any $n \geq n_0$, any
tournament $G$ on at least $2(1+\alpha)n$ vertices contains any directed tree
$T$ on $n$ vertices.
\item Let $\Delta$ be any positive integer. Then there exists $n_0 \in
\mathbb{N}$ such that for any $n \geq n_0$, any tournament $G$ on at least
$(1+\alpha)n$
vertices contains any directed tree $T$ on $n$ vertices with $\Delta(T) \leq
\Delta$.
\end{enumerate}}

\medskip
The proofs of each of the two statements of the theorem are very similar, so to
avoid repetition we shall prove the first statement, explaining in footnotes
where the proof of the second statement differs.

\subsection{Partitioning the vertices of $G$.}

As in the last section, we consider the set $A$ of all positive values of
$\alpha$ such
that the theorem holds. So $\alpha' \in A$ if and only if there exists $n_0$ such
that for any $n \geq n_0$, any tournament on at least $2(1+\alpha')n$ vertices
contains
any tree on $n$ vertices. So if $\alpha' \in A$ and $\alpha'' >
\alpha'$ then $\alpha'' \in A$, and also $1/2 \in A$
by Theorem~\ref{best_so_far}.
Thus we may define $a_{\textrm{inf}} = \inf A$, and then the
theorem is equivalent to the statement that $a_{\textrm{inf}} = 0$. So suppose $a_{\textrm{inf}} >
0$, and choose constants $$1/n_0 \ll 1/n'_0 \ll \mu \ll \nu \ll \eta \ll
1/\Delta' \ll \gamma
\ll a_{\textrm{inf}}.$$ Let $\alpha = a_{\textrm{inf}} - \mu$, so $\alpha \leq 1/2$, and we may
assume that $\gamma \ll \alpha$. Then $\alpha + 2 \mu \in A$, and so
for any $n' \geq n'_0$, any tournament on at least $2(1 +\alpha + 2\mu)n'$
vertices
contains any tree on $n'$ vertices. We shall prove that if $n \geq n_0$, any
tournament $G$ on at least $2(1+\alpha) n$ vertices contains any tree on $n$
vertices. This proves that $\alpha \in A$,
giving a contradiction to our assumption that $a_{\textrm{inf}} > 0$, and so proving the
theorem.\footnote{For the bounded degree case, fix any value of $\Delta$, and
here $A = A(\Delta)$ is defined by $\alpha' \in A$ if and only if there exists $n_0$ such
that for any $n \geq n_0$, any tournament on at least $(1+\alpha')n$ vertices
contains any tree~$T$ on~$n$ vertices with $\Delta(T) \leq \Delta$.
So if $\alpha' \in A$ and $\alpha'' > \alpha'$ then $\alpha'' \in A$, and also
$2 \in A$ by Theorem~\ref{best_so_far}. Thus we may define $a_{\textrm{inf}} = \inf A$;
then the theorem is equivalent to the statement that $a_{\textrm{inf}} = 0$. So suppose
$a_{\textrm{inf}} > 0$, and choose constants $1/n_0 \ll 1/n'_0 \ll \mu \ll \nu \ll \eta
\ll 1/\Delta' \ll  \gamma \ll 1/\Delta,a_{\textrm{inf}}$. Let $\alpha = a_{\textrm{inf}} - \mu$,
so $\alpha < 2$, and we may assume that $\gamma \ll \alpha$. Then $\alpha + 2\mu \in A$, so for any $n'
\geq n'_0$, any
tournament on at least $(1 +\alpha + 2\mu)n'$ vertices contains any tree $T$ on
$n'$
vertices with $\Delta(T) \leq \Delta$. Using this, we shall prove that if $n
\geq n_0$, any tournament $G$ on at least $(1+\alpha) n$ vertices contains any
tree $T$ on $n$ vertices with $\Delta(T) \leq \Delta$. This proves that $\alpha
\in A$, giving a contradiction to our assumption that $a_{\textrm{inf}} > 0$, and so
proving the theorem.}

So let $G$ be a tournament on at least $2(1+\alpha)n$ vertices\footnote{For the
bounded degree case, instead let $G$ be a tournament on at least $(1+\alpha) n$
vertices.}. If $|G| \geq 3n$ then by Theorem~\ref{best_so_far},~$G$ contains any
directed tree $T$ on $n$ vertices. So we may assume $|G| < 3n$. We shall use an
algorithm which keeps track of an ordered family
$\mathcal{S}^\tau$ of disjoint subsets of $V(G)$, and a set~$B^\tau$ of bad
edges of $G$, at each time $\tau$. Initially, let $\mathcal{S}^1 = (V(G))$, and
let $B^1 =
\emptyset$. Then at time $\tau \geq 1$, we have $\mathcal{S}^\tau = (S_1^\tau,
\dots, S_\tau^\tau)$, and the algorithm proceeds as follows.
\begin{enumerate}
\item Let $S_\ell^\tau$ be a largest member of $\mathcal{S^\tau}$. If
$|S_\ell^\tau| < \gamma n$, then terminate.
\item If $G[S_\ell^\tau]$ is a robust $(\mu, \nu)$-outexpander with
$\delta^0(G[S_\ell^\tau]) \geq \eta n$, then terminate.
\item If some $v \in S_\ell^\tau$ has $d^+_{G[S_\ell^\tau]}(v) < \eta n$, then let
$$\mathcal{S}^{\tau+1} = (S_1^\tau, \dots, S_{\ell-1}^\tau, S_\ell^\tau \sm
\{v\},
\{v\}, S_{\ell+1}^\tau, \dots, S_{\tau}^\tau), $$
let $B^{\tau+1} = B^\tau \cup E(\{v\} \rightarrow S_\ell^\tau \sm \{v\})$, and
proceed to step (6).
\item Similarly, if some $v \in S_\ell^\tau$ has $d^-_{G[S_\ell^\tau]}(v) < \eta
n$, then let
$$\mathcal{S}^{\tau+1} = (S_1^\tau, \dots, S_{\ell-1}^\tau, \{v\}, S_\ell^\tau \sm
\{v\}, S_{\ell+1}^\tau, \dots, S_{\tau}^\tau),$$
let $B^{\tau+1} = B^\tau \cup E(S_\ell^\tau \sm \{v\} \rightarrow \{v\})$, and
proceed to step (6).
\item If $G[S_\ell^\tau]$ is not a robust $(\mu, \nu)$-outexpander then apply
Lemma~\ref{not_robust_expander_split} to partition the vertices of $S_\ell^\tau$
into sets $S'$ and $S''$ such
that $\nu|S_\ell^\tau| \leq |S'|, |S''| \leq (1-\nu)|S_\ell^\tau|$ and at most
$4\mu |S_\ell^\tau|^2$ edges of $G[S_\ell^\tau]$ are directed from $S''$ to $S'$.
Then let
$$\mathcal{S}^{\tau+1} = (S_1^\tau, \dots, S_{\ell-1}^\tau, S', S'',
S_{\ell+1}^\tau,
\dots, S_{\tau}^\tau)$$
and let $B^{\tau+1} = B^\tau \cup E(S'' \rightarrow S')$.
\item Finally, for each $i \in [\tau+1]$, delete from $S_i^{\tau+1}$ any vertex
$v$ which lies in more than $\sqrt{\eta}n$ edges of $B^{\tau+1}$.
\end{enumerate}

At any step, if the algorithm does not terminate at step (1) or (2), then the
condition of one of steps (3), (4) and (5) must hold. Therefore at each time
$\tau$, either the algorithm terminates or $|\mathcal{S}^\tau|$ increases from
$\tau$ to $\tau+1$ (in forming $\mathcal{S}^{\tau+1}$) by reducing the size of
the largest piece. Therefore the algorithm must terminate\COMMENT{COMMENT: Is
this obvious? Formally, let $D$ be the set of vertices which are deleted at any
point in the algorithm. Then at time $\tau$, we have $\tau$ members of
$\mathcal{S}^\tau$. The only way any one of these members can have become empty
is by having all of its elements deleted, and so at most $|D|$ members of
$\mathcal{S^\tau}$ can be empty. But then there are $|G| - |D|$ undeleted
vertices split between $\tau - |D|$ non-empty disjoint sets in
$\mathcal{S}^\tau$, so we must have $|G| \geq \tau$.} at some time $\tau_{\textrm{end}}
\leq |G| \leq 3n$.

Now, at any time $\tau$ at which the algorithm does not terminate, the algorithm
will split the set $S^\tau_\ell$ in precisely one of steps (3), (4) and (5). We
next show that the split in step (5) will occur for at most $3/\gamma \nu$ times
$\tau < \tau_{\textrm{end}}$. This is because any set obtained by a split in step (5)
must have size at least $\gamma \nu n$ (since $|S^\tau_\ell| \geq \gamma n$, and
the sets $S', S''$ obtained have $|S'|, |S''| \geq \nu |S^\tau_\ell|$), and so
at most $|G|/\gamma \nu n\leq 3/\gamma \nu$ such sets can be
obtained.

Next, we show that when the algorithm terminates at time $\tau_{\textrm{end}}$, most
vertices lie in one of the sets $S_i^\tau$, or equivalently that only a few vertices
have been deleted. To do this, note that at each time $\tau \leq \tau_{\textrm{end}}$,
the number of edges added to form $B^{\tau+1}$ from $B^\tau$ is at most $\eta n$
if the algorithm carried out the split in step (3) or (4), and at most $4\mu
|G|^2 \leq 36\mu n^2$ if the algorithm carried out the split in step (5). Since
$\tau_{\textrm{end}} \leq 3n$, and the split in step (5) is carried out in at most
$3/\gamma \nu$ of these steps, the number of bad edges at time $\tau_{\textrm{end}}$ is
at most $3 \eta n^2 + 108
\mu n^2/\nu \gamma \leq 4 \eta n^2$. Since $B^1 \subseteq \dots \subseteq
B^{\tau_{\textrm{end}}}$, any vertex of $G$ which was ever deleted in step (6) must lie
in at least $\sqrt{\eta} n$ edges of $B^{\tau_{\textrm{end}}}$, and so at most $8
\sqrt{\eta} n$ vertices of $G$ can have been deleted in step (6) over the entire
course of the algorithm. Let $G'$ be the restriction of~$G$ to the undeleted
vertices at time $\tau_{\textrm{end}}$, so $G' = G[\bigcup \mathcal{S}^{\tau_{\textrm{end}}}]$.
Then $G'$ is a tournament and
\begin{equation} \label{eq:size_of_G'}
|G'| \geq |G|- 8 \sqrt{\eta} n.
\end{equation}

Our approach now depends on whether the algorithm terminated in step (1) or (2).
If the algorithm terminated in step (1), then for each $i
\in [{\tau_{\textrm{end}}}]$ we have $|S_i^{\tau_{\textrm{end}}}| < \gamma n$. We shall show that
in this case $G'$ is $2\gamma$-almost-transitive. Indeed, order the vertices of
$G'$ as $v_1,
v_2, \dots, v_{|G'|}$ in the same order as in $\mathcal{S}^{\tau_{\textrm{end}}}$, i.e.
beginning with
all the vertices of $S^{\tau_{\textrm{end}}}_1$, then the vertices of $S^{\tau_{\textrm{end}}}_2$,
and so forth.
Then any edge $v_j \rightarrow v_i$ where $j>i$ either lies in \COMMENT{COMMENT:
By induction, at any time $\tau$ any edge from $S_j^\tau$ to $S_i^\tau$ where
$j>i$ is in $B^\tau$, since in forming $\mathcal{S}^{\tau+1}$ the new edges
formed of this type are those from $\{v\}$ to $S_\ell^\tau \sm \{v\}$ if we
split in step (3), those from $S_\ell^\tau \sm \{v\}$ to $\{v\}$ if we split in
step(4), and those from $S''$ to $S'$ if we split in step (5). In any case,
these are the edges that we add to $B^\tau$ to form $B^{\tau+1}$.}
$B^{\tau_{\textrm{end}}}$ or has both endvertices in the same $S_i^{\tau_{\textrm{end}}}$. So the
total number of such edges is at most
$$4 \eta n^2 + \sum_{S \in \mathcal{S}^{\tau_{\textrm{end}}}}
\binom{|S|}{2} \leq 4 \eta n^2 + \sum_{S \in \mathcal{S}^{\tau_{\textrm{end}}}}
\frac{\gamma n|S|}{2} \leq 4 \eta n^2 + \frac{3\gamma n^2}{2} \leq 2\gamma n^2.$$
Since in both the unbounded degree case and the bounded degree case we
have
\begin{equation*}
|G'| \geq (1+\alpha/2)n,
\end{equation*} by (\ref{eq:size_of_G'}), $G'$ is indeed
$2\gamma$-almost-transitive, and by Lemma~\ref{almost_trans} $G'$
contains a copy of $T$, which is also a copy of $T$ in $G$.

\subsection{Partitioning the vertices of $T$.}

We may therefore assume that the algorithm terminated in step (2) at some time
${\tau_{\textrm{end}}}$; when for some $S^{\tau_{\textrm{end}}}_i$ with $|S^{\tau_{\textrm{end}}}_i| \geq
\gamma n$, $G[S^{\tau_{\textrm{end}}}_i]$
is a $(\mu, \nu)$-robust outexpander with $\delta^0(G[S^{\tau_{\textrm{end}}}_i]) \geq \eta
n$. For this $i$, let
$S = S^{\tau_{\textrm{end}}}_i$, let $S^+ = \bigcup_{i<j\leq {\tau_{\textrm{end}}}}
S^{\tau_{\textrm{end}}}_j$ and let $S^- =
\bigcup_{1\leq j<i} S^{\tau_{\textrm{end}}}_j$. Then $|S^+ \cup S^- \cup S| = |G'|$.
Also, if $u \in S^+$ and $v \in S \cup S^-$ then $u \in S^{\tau_{\textrm{end}}}_j$, $v \in
S^{\tau_{\textrm{end}}}_\ell$ for some $j > \ell$, and so if $u \rightarrow v$ then this
edge is in $B^{\tau_{\textrm{end}}}$. So any vertex $u \in S^+$ has at most
$\sqrt{\eta}n$ outneighbours in $S \cup S^-$, since $u$ was not deleted at any
stage of the algorithm. Similarly each vertex of $S$ has at most $\sqrt{\eta}n$
outneighbours in $S^-$ and inneighbours in $S^+$, and each vertex of $S^-$ has
at most $\sqrt{\eta}n$ inneighbours in $S^+ \cup S$. Define $\beta, \beta^+,
\beta^-$ by $|S| = \beta |G'| $, $|S^+| = \beta^+ |G'|$, and $|S^-|
=\beta^-|G'|$, so $\beta + \beta^+ + \beta^- = 1$ and $\beta \geq \gamma n/|G'|
\geq \gamma /3$.

Suppose first that $\beta^+$ and $\beta^-$ are both small. More precisely,
$\beta^+, \beta^- \leq \alpha\beta^2/20$, and so $\beta \geq 1-\alpha/10$. Then we
shall find a copy of $T$ in $G[S]$ (and therefore in $G$). Indeed, $T$ is a tree
on $n$ vertices, and $G[S]$ is a $(\mu, \nu)$-robust outexpander with
$\delta^0(G[S]) \geq \eta n$. Furthermore, $$|S|
= \beta |G'| \stackrel{(\ref{eq:size_of_G'})}{\geq} (2+2\alpha-8\sqrt{\eta})\beta n  \geq
(2+\alpha)(1-\frac{\alpha}{10})n \geq 2(1+\frac{\alpha}{4}) n$$ and so by
Lemma~\ref{robust_expander_case} $G[S]$ (and therefore $G$) contains a copy of
$T$.\footnote{For the bounded degree case, $|S| = \beta |G'|
\stackrel{(\ref{eq:size_of_G'})}{\geq}
(1+\alpha-8\sqrt{\eta})(1-\alpha/10)n \geq (1+\alpha/4) n$, and so $G[S]$ (and
therefore $G$) contains a copy of $T$ by
Lemma~\ref{robust_expander_case_bounded_deg}.}

So we may assume that at least one of $\beta^+$ and $\beta^-$ is greater than
$\alpha\beta^2/20$, so in particular, $\beta \leq 1-\alpha\beta^2/20$. We next
split the vertices of $T$ according to the values of
$\beta^+$ and $\beta^-$.

\emph{Case 1: $\beta^-$ is large but $\beta^+$ is small.} More precisely,
$\beta^+ \leq \alpha\beta^2/20$ and $\beta^- > \alpha\beta^2/20$.
Then we partition the vertex set of $T$ into $T^-$ and $T^0$, where every edge
of $T$ between $T^-$ and $T^0$ is directed from $T^-$ to $T^0$, and $|T^-| =
\beta^-(1-\alpha\beta)n$. We can form $T^0$ greedily by successively removing a
sink vertex from $T$ and adding it to $T^0$.
Since $\beta^+ + \beta + \beta^- = 1$,
\begin{equation*} \label{eq:sec5_T'case1}
|T^0| = n - |T^-| = \beta n (1+\alpha -
\alpha\beta) + (1-\alpha\beta)\beta^+n \leq \beta n(1+\alpha - \alpha\beta) +
\alpha\beta^2 n/20.
\end{equation*}

\emph{Case 2: $\beta^+$ is large but $\beta^-$ is small.} More precisely,
$\beta^- \leq \alpha\beta^2/20$ and $\beta^+ > \alpha\beta^2/20$.
Then we similarly partition the vertex set of $T$ into~$T^0$ and $T^+$,
where every edge of
$T$ between~$T^0$ and $T^+$ is directed from $T^0$ to $T^+$, and $|T^+| =
\beta^+(1-\alpha\beta)n$. Again $|T^0| = n - |T^+| \leq \beta n (1+\alpha -
\alpha\beta)
+\alpha\beta^2 n/20 $.

\emph{Case 3: $\beta^+$ and $\beta^-$ are both large.} More precisely, $\beta^+,
\beta^- > \alpha\beta^2/20$. Then we partition\COMMENT{COMMENT: This is possible by
choosing first $T^-$ then $T^+$ as in the previous cases.} the vertex set of $T$
into pieces $T^-$, $T^0$ and $T^+$ such that all edges of $T$ between $T^-$ and
$T^0$ are directed from $T^-$ to $T^0$, all edges of $T$ between $T^0$ and $T^+$
are directed from $T^0$ to $T^+$ and all edges of $T$ between $T^-$ and $T^+$ are
directed from $T^-$ to $T^+$. Also $|T^+| = \beta^+(1-\alpha\beta)n$ and $|T^-|
= \beta^-(1-\alpha\beta)n$, so $|T^0| = \beta(1+\alpha-\alpha\beta)n$.

Note that in each of the three cases $T^0$ satisfies
$|T^0| \geq \beta(1+\alpha-\alpha\beta)n$ and
\begin{equation} \label{eq:sec5_T'case3}
|T^0| \leq \beta  (1+\alpha - \alpha\beta)n
+\alpha\beta^2 n/20 \le \beta (1+\alpha)n-\frac{\alpha \beta^2 n}{2}.
\end{equation}

\subsection{Embedding $T$ in $G$.}

Having partitioned the vertices of $G'$ into three sets $S, S^+$ and $S^-$, and
the tree $T$ into three forests $T^+$, $T^0$, $T^-$, we now complete the proof by
embedding $T$ in $G$, with $T^-, T^0$ and $T^+$ embedded in $G[S^-], G[S]$ and
$G[S^+]$ respectively. Indeed, the fact that $G[S]$ is a robust $(\mu,
\nu)$-outexpander will enable us to embed slightly more vertices in $G[S]$ than
the $\beta n$ that would be embedded in $G[S]$ if the vertices of $T$ were
distributed proportionately amongst $G[S]$, $G[S^+]$ and $G[S^-]$. This gives us
some leeway for embedding $T^+$ and $T^-$ in $G[S^+]$ and $G[S^-]$ respectively,
which by our choice of $\alpha$ is sufficient to successfully complete these
embeddings.

So let $T_1^-, \dots, T_x^-$ be the component
subtrees of $T^-$, let $T_1^+, \dots, T_y^+$ be the component subtrees of $T^+$,
and let $T_1, \dots, T_z$ be the component subtrees of $T^0$. Let the
\emph{contracted tree}~$T_{\textrm{con}}$ be formed from $T$ by contracting each $T_i^+,
T_i^-$ and $T_i$ to a single vertex.

To begin the embedding, we embed into $G[S]$ every $T_i$ satisfying $|T_i| \geq
n/\Delta'$. Note that there are at most $\Delta'$ such $T_i$. Also, the union of
all such $T_i$ is a forest on at most $|T^0|$ vertices, and the tournament
$G[S]$ is a robust $(\mu, \nu)$-outexpander on
$$\beta |G'|  \stackrel{(\ref{eq:size_of_G'})}{\geq}
\beta (2+2\alpha - 8\sqrt{\eta})n
\stackrel{(\ref{eq:sec5_T'case3})}{\geq} 2\left(1+\frac{\alpha\beta}{10}\right)
|T^0| \geq 2(1+\gamma^2)|T^0|$$
vertices with $\delta^0(G[S]) \geq
\eta n$, and hence
 $G[S]$ contains a copy of
this forest by Lemma~\ref{robust_expander_case}.\footnote{For the bounded degree
case, $|S| \geq (1+\gamma^2) |T^0|$ by a similar calculation, and so $G[S]$
contains a copy of this forest by Lemma~\ref{robust_expander_case_bounded_deg}.}

Now, choose an order of the vertices of $T_{\textrm{con}}$, beginning with the at most
$\Delta'$ vertices corresponding to the $T_i$ which we have just embedded, and
such that any vertex of $T_{\textrm{con}}$ has at most $\Delta'$ neighbours preceding it
in this order. (To do this, choose one of the $\Delta'$ vertices corresponding
to the $T_i$ which have already been embedded, and then choose any ancestral
ordering of the vertices of $T_{\textrm{con}}$, beginning with the chosen vertex,
so every vertex has at most one
neighbour preceding it in this order.
Now move the remaining $\Delta'-1$ vertices corresponding to the $T_i$ which
have already been embedded to the
front of this order; then every vertex gains at most $\Delta'-1$ preceding
neighbours.)
We shall proceed through the remaining vertices of $T_{\textrm{con}}$ in this order, at
each step embedding the tree $T_i, T_i^+$ or $T_i^-$ corresponding to the
current vertex of $T_{\textrm{con}}$ in the unoccupied vertices of the tournament $G[S],
G[S^+]$ or $G[S^-]$ respectively.

So suppose first that the current vertex $t^*$ of $T_{\textrm{con}}$ corresponds to some
$T_i$. Since $T_i$ has not already been embedded, we know that $|T_i| \leq
n/\Delta'$. Also, since $t^*$ has at most $\Delta'$ neighbours preceding it in
$T_{\textrm{con}}$, the vertices of $T_i$ have at most $\Delta'$ neighbours outside $T_i$
which have already been embedded. Since $T_i$ is a component of $T^0$, each
of these neighbours of vertices in $T_i$ lies either in $T^-$ (in which
case it is an
inneighbour) or in $T^+$ (in which case it is an outneighbour). So let $t_1^-,
\dots, t_p^-$ be the vertices in $T^-$ which are inneighbours of some vertex in
$T_i$ and which have previously been embedded, and let $v_1^-, \dots, v_p^-$
be the vertices of $G'[S^-]$ to which
$t_1^-, \dots, t_p^-$ were embedded. Similarly, let $t_1^+, \dots, t_q^+$
be the vertices in $T^+$ which are outneighbours of some vertex in $T_i$ and
which have previously been embedded, and
let $v_1^+, \dots, v_q^+$ be the vertices of $G'[S^+]$ to which $t_1^+, \dots,
t_q^+$ were embedded. Finally let $S^*$ be the set of unoccupied vertices
in $S \cap N^+(v_1^-, \dots, v_p^-) \cap N^-(v_1^+, \dots, v_q^+)$. Then we wish
to embed~$T_i$ in $S^*$. For this, note that
\begin{align*}
|S^*| &\geq |S| - (p+q)\sqrt{\eta}n - |T^0|
 \stackrel{(\ref{eq:sec5_T'case3})}{\geq} \beta|G'| - \Delta' \sqrt{\eta}n  -
(\beta (1+\alpha)n - \alpha\beta^2n/2)
\\ &\stackrel{(\ref{eq:size_of_G'})}{\geq} \beta n(1+\alpha ) - (8+\Delta')
\sqrt{\eta}n- \beta (1+\alpha )n+
\alpha\beta^2 n/2
 \geq \alpha\beta^2 n/3 \geq 3n/\Delta' \geq 3|T_i|.
\end{align*}
Note that this calculation is valid for both the bounded degree
case and the unbounded degree case, with plenty of room to spare in the
unbounded case.
So by Theorem~\ref{best_so_far}, $G[S^*]$ contains a copy of $T_i$, to which
we embed $T_i$.

Alternatively, if the current vertex of $T_{\textrm{con}}$ corresponds to some $T_i^-$,
then similarly the vertices of $T_i^-$ have at most $\Delta'$ neighbours outside
$T_i^-$ which have already been embedded, all of which are outneighbours. As
before we let $v_1, \dots, v_r$ be
the vertices of $G'[S \cup S^+]$ to which these vertices have been embedded, and
let $S^*$ be
the set of unoccupied vertices of $S^- \cap N^-(v_1, \dots, v_r)$. Note that at
most $|T^-| - |T_i^-|$ vertices of $T^-$ have already been embedded. Since some
$T_i^-$ exists we have
\begin{align} \label{eq:final}
|S^*| &\geq |S^-| - r\sqrt{\eta}n - (|T^-| - |T_i^-|)
\stackrel{(\ref{eq:size_of_G'})}{\geq} \beta^-(2+2\alpha)n -
(8+\Delta') \sqrt{\eta}n -
\beta^-(1-\alpha\beta)n
+ |T_i^-|
\nonumber\\ &\geq \beta^-(1+2\alpha+\alpha\beta/2)n + |T_i^-|.
\end{align}
In the final line we used the fact that $\beta^- \geq \alpha\beta^2/20$ and
$\beta \geq \gamma/3$ (so $\eta, 1/\Delta' \ll \gamma, \beta, \beta^-$).
So $|S^*| \geq
2(1+\alpha+2\mu)|T_i^-|$.
Therefore if $|T_i^-| \geq \beta^- n/2$, then $|T_i^-| \geq \alpha\beta^2n/40
\geq \alpha\gamma^2n/360 \geq n_0'$, and so we can embed $T_i^-$ in $G[S^*]$
by our choice of $n_0'$. On the other
hand, if $|T_i^-| < \beta^- n/2$ then $|S^*| \geq 3|T_i^-|$ by~(\ref{eq:final}),
and so we can embed
$T_i^-$ in $G[S^*]$ by Theorem~\ref{best_so_far}.\footnote{For the
bounded degree case
\begin{align*}
|S^*| &\geq |S^-| - r\sqrt{\eta}n - (|T^-| - |T_i^-|)
\stackrel{(\ref{eq:size_of_G'})}{\geq} \beta^-(1+\alpha)n -
(8+\Delta') \sqrt{\eta}n -
\beta^-(1-\alpha\beta)n +
|T_i^-|
\\ &\geq \beta^-(\alpha+\alpha\beta/2)n + |T_i^-|.
\end{align*}
So $|S^*| \geq
(1+\alpha+2\mu)|T_i^-|$.
Therefore if $|T_i^-| \geq \beta^-\alpha n/2$, then
$|T_i^-| \geq n_0'$, and so we can embed $T_i^-$ in $G[S^*]$
by our choice of $n_0'$. On the other
hand, if $|T_i^-| < \beta^-\alpha n/2$ then $|S^*| \geq 3|T_i^-|$,
and so we can embed $T_i^-$ in $G[S^*]$ by Theorem~\ref{best_so_far}.
}

Finally, if the current vertex of $T_{\textrm{con}}$ corresponds to some $T_i^+$, we
embed $T_i^+$ in the unoccupied vertices of $S^+$ by a similar method to the
method used to embed some $T_i^-$ in the
unoccupied vertices of $G[S^-]$. We continue in this
manner until we have embedded the $T_i$,~$T_i^+$ or $T_i^-$ corresponding to
each vertex of $T_{\textrm{con}}$, at which point we will have obtained an embedding of
$T$ in $G$, completing the proof. At each stage in this proof we had `room to
spare' in our choices, and so the fact that the expressions for $|T_i|$,
$|T_i^+|$ and $|T_i^-|$ and other such expressions may not be integers is not a
problem.
\endproof

\section*{Acknowledgements}

We would like to thank the anonymous referees for their many helpful suggestions and comments.

\medskip

\noindent
{\footnotesize
\noindent
Daniela K\"uhn, Deryk Osthus, \\School of Mathematics,
\\University of Birmingham, \\Birmingham, \\B15 2TT, \\United Kingdom, \\
\{{\tt kuehn,osthus}\}{\tt @maths.bham.ac.uk } }

\medskip
\noindent
{\footnotesize
\noindent
Richard Mycroft, \\ School of Mathematical Sciences,
\\ Queen Mary, \\ University of London, \\London, \\E1 4NS, \\United Kingdom, \\
{\tt r.mycroft@qmul.ac.uk } }
\end{document}